\documentclass[ijoc, nonblindrev]{informs3} 
\OneAndAHalfSpacedXII 

\usepackage{diagbox}
\usepackage{booktabs}
\usepackage[ruled, vlined, linesnumbered]{algorithm2e}
\usepackage{setspace}
\usepackage{makecell}
\usepackage{multirow}
\usepackage{xcolor}
 \usepackage{mathrsfs}
\def\myright#1{\unskip\nobreak\hfill#1}

\SetCommentSty{mycommfont}

\usepackage{natbib}
 \bibpunct[, ]{(}{)}{,}{a}{}{,}%
 \def\bibfont{\small}%

 \usepackage{mathtools}
\DeclarePairedDelimiter{\ceil}{\lceil}{\rceil}

\TheoremsNumberedThrough     

\EquationsNumberedThrough    

\MANUSCRIPTNO{JOC-2023-04-OA-0128} 

\begin{document}


\RUNAUTHOR{Zhang et al.} 

\RUNTITLE{DIDP for SUALBP}

\TITLE{Domain-Independent Dynamic Programming and Constraint Programming Approaches for \\ Assembly Line Balancing Problems with Setups}

\ARTICLEAUTHORS{%
\AUTHOR{Jiachen Zhang, J. Christopher Beck}
\AFF{Department of Mechanical and Industrial Engineering, University of Toronto, Toronto, Ontario M5S 3G8, Canada, \\ \EMAIL{\{jasonzjc,  jcb\}@mie.utoronto.ca} \URL{}}
} 

\ABSTRACT{%
We propose domain-independent dynamic programming (DIDP) and constraint programming (CP) models to exactly solve type-1 and type-2 assembly line balancing problem with sequence-dependent setup times (SUALBP). The goal is to assign tasks to assembly stations and to sequence these tasks within each station, while satisfying precedence relations specified between a subset of task pairs. Each task has a given processing time and a setup time dependent on the previous task on the station to which the task is assigned. The sum of the processing and setup times of tasks assigned to each station constitute the station time and the maximum station time is called the cycle time. For type-1 SUALBP, the objective is to minimize the number of stations, given a maximum cycle time. For type-2 SUALBP, the objective is to minimize the cycle time, given the number of stations. On a set of diverse SUALBP instances, experimental results show that our approaches significantly outperform the state-of-the-art mixed integer programming models for SUALBP-1. For SUALBP-2, the DIDP model outperforms the state-of-the-art exact approach based on logic-based Benders decomposition. By closing 76 open instances for SUALBP-2, our results demonstrate the promise of DIDP for solving complex planning and scheduling problems.



}%


\KEYWORDS{domain-independent dynamic programming, constraint programming, assembly line balancing}

\maketitle

%




\section{Introduction}

The \emph{simple assembly line balancing problem} (SALBP)
is a well-known production planning problem  \citep{becker2006survey}, with many applications in the production of automotive and household items \citep{boysen2022assembly}.
The SALBP is often the last step of production, performing the final assembly of a product from previously manufactured parts. Commonly used for high-volume production of standardized products, an assembly line is a series of workstations connected by a transport system such as a conveyor and the production process is centered on work pieces that move through the workstations in a fixed order. At each station, a set of tasks is performed sequentially. Due to the technical structure of the products, the set of tasks must adhere to a partial precedence ordering: if $(a, b)$ represents a precedence relation between task $a$ and $b$, then $b$ must be assigned to the same or a later station than $a$. If assigned to the same station, $b$ must be after $a$ in the task sequence. In a feasible production plan, each task is assigned to exactly one station and the precedence constraints are satisfied. The set of tasks assigned to each station is the station workload and the sum of task times is the station time. The maximum station time among all stations is called the cycle time.

In many real production lines, setup operations such as tool changes,
curing, or cooling processes are required between consecutive tasks \citep{kumar2013assembly}. The \emph{assembly line balancing problem with setups} (SUALBP) incorporates setup times in the SALBP \citep{andres2008balancing}. Specifically, for all the tasks assigned to the same station, sequence-dependent setup times need to be considered in the task schedule. Scholl
et al. \citeyearpar{scholl2013assembly} modify the problem of Andrés et al. \citeyearpar{andres2008balancing} and distinguish forward and backward setup times. A forward setup time occurs when a
task directly follows another task on the same station
and cycle. A backward setup time is performed
between the last task in a cycle and the first task in the
following cycle on the same station \citep{zohali2022solving}.

In the literature, SUALBP has also been referred to as \emph{general assembly line balancing problem with setups} \citep{martino2010heuristic} and \emph{setup assembly line balancing and scheduling problem} \citep{scholl2013assembly}. There are two main types of SUALBP: SUALBP-1 minimizes the number of stations given a fixed cycle time and SUALBP-2 minimizes the cycle time given a fixed number of stations. Due to different application scenarios, other variants also exist in the literature \citep{wee1982assembly, becker2006survey}.

There has been a variety of approaches applied to the SUALBP-1 variants including \emph{mixed-integer programming} (MIP) \citep{akpinar2014modelinga}, constraint programming \citep{guner2023constraint}, hybrid bee colony algorithms \citep{akpinar2014modelingb}, Benders decomposition \citep{akpinar2017combinatorial}, simulated annealing  \citep{ozcan2019balancing}, and variable
neighborhood search \citep{yang2020modelling}. For SUALBP-2 variants, MIP models \citep{tang2016balancing}, genetic algorithms combined with dynamic programming \citep{yolmeh2012efficient}, and meta-heuristics \citep{csahin2017increasing} have been developed.

For the standard SUALBP-1 and SUALBP-2, Andres et al. \citeyearpar{andres2008balancing} develop a mathematical model and multiple heuristics to solve SUALBP-1.
Martino and Pastor 
\citeyearpar{martino2010heuristic} develop heuristic algorithms for the SUALBP-1. Later, Seyed-Alagheband et al. \citeyearpar{seyed2011simulated} develop a MIP model based
on the formulation of Andrés et al. \citeyearpar{andres2008balancing} and propose a simulated annealing meta-heuristic to solve SUALBP-2.  Scholl et al. \citeyearpar{scholl2013assembly} propose a new MIP model and multiple heuristics to solve a SUALBP-1 variant.
Esmaeilbeigi et al. \citeyearpar{esmaeilbeigi2016new}
propose improved MIP models for both SUALBP-1 and SUALBP-2 that represent the state-of-the-art exact techniques for the former. More recently, Zohali et al. \citeyearpar{zohali2022solving} develop the state-of-the-art exact approach to SUALBP-2 using logic-based Benders decomposition.

In this work, we develop novel \emph{domain-independent dynamic programming} (DIDP) and \emph{constraint programming} (CP) models for both types of SUALBP. The superior performance of the proposed DIDP models against the state-of-the-art approaches demonstrates the promise of this emerging exact method for solving complex planning and scheduling problems.

\paragraph{Main contributions.} 
Our contributions are four novel optimization models (DIDP and CP) for SUALBP-1 and SUALBP-2, and new state-of-the-art results of the proposed DIDP models compared to the best exact algorithms in the literature. We successfully close 76 open instances for SUALBP-2 and detect potential flaws in empirical results of Zohali et al. \citeyearpar{zohali2022solving}. We also investigate a local improvement algorithm for the DIDP model of SUALBP-2.

\paragraph{Outline of the paper.} In Section 2, we define the problem, summarize the notation, present the state-of-the-art MIP models and exact algorithms, and introduce domain-independent dynamic programming and constraint programming. The proposed DIDP and CP models for SUALBP-1 are presented in Section 3. Similarly, the DIDP and CP models for SUALBP-2 are presented in Section 4. A local improvement algorithm to solve the DIDP model of SUALBP-2 is briefly introduced in Section 5. In Section 6, experimental results on benchmark instances are presented, followed by discussions regarding potential flaws in previous work. Finally, we conclude this paper in Section 6.

\section{Background}

In this section, we define the assembly line balancing problem with setups, summarize the notation, present the state-of-the-art MIP models, briefly introduced the state-of-the-art logic-based Benders decomposition for SUALBP-2, and review domain-independent dynamic programming and constraint programming methodologies.

\subsection{Problem Definition and Notation}

\begin{table}[htbp]
    \centering
    \caption{Notation and definition for SUALBP-1 and SUALBP-2 \citep{esmaeilbeigi2016new}.}
    \label{table:def}
        \begin{tabular}{p{.15\textwidth}p{.78\textwidth}}
        \hline
        Notation & Definition \\ \hline
        $n$ & the number of tasks \\
        $V$ & set of tasks, i.e., $V = \{1,2,...,n\}$ \\
        $i,j,v$ & index for the tasks \\
        $k$ & index for the stations \\
        $\mathcal{E}$ & set of all precedence relations \\
        $t_{i}$ & execution time for task $i\in V$ \\
        $P_{i}(P_{i}^{*})$ & set of direct (all) predecessors of task $i\in V$ \\
        $F_{i}(F_{i}^{*})$ & set of direct (all) followers of task $i\in V$ \\
        $\overline{c}(\underline{c})$ & upper (lower) bound on the cycle time \\
        $\overline{m}(\underline{m})$ & upper (lower) bound on the number of stations \\
        $E_{i}$ & earliest station for task $i\in V$, e.g., $E_{i} = \ceil{\frac{t_{i} + \sum_{j\in P_{i}^{*}}t_{j}}{\overline{c}}}$ \\
        $L_{i}$ & latest station for task $i\in V$, e.g., $L_{i} = \overline{m} + 1 - \ceil{\frac{t_{i} + \sum_{j\in F_{i}^{*}}t_{j}}{\overline{c}}}$ \\
        $KD$($KP$) & set of definite (possible) stations, i.e., $KD = \{1,...,\underline{m}\}$, and $KP = \{\underline{m}+1, ..., \overline{m}\}$ \\
        $K$ & set of all stations, i.e., $K = KD \cup KP$ \\
        $FS_{i}$ & set of stations to which task $i\in V$ can be assigned, i.e., $FS_{i} = \{E_{i}, E_{i}+1, ..., L_{i}\}$ \\
        $FT_{k}$ & set of tasks which can be assigned to station $k\in K$, i.e., $FT_{k} = \{i\in V | k\in FS_{i}\}$ \\
        $A_{i}$ & set of tasks that cannot be assigned to the station to which task $i$ is assigned, e.g., $A_{i} = \{j\in V | FS_{j} \cap FS_{i} = \emptyset\}$ \\
        $F_{i}^{F}(P_{i}^{F})$ & set of tasks which may directly follow (precede) task $i$ in forward direction, i.e., $F_{i}^{F} = \{j \in V-(F_{i}^{*} - F_{i}) - P_{i}^{*} - A_{i} - \{i\}\}$ and $P_{i}^{F} = \{j\in V|i\in F_{j}^{F}\}$ \\
        $F_{i}^{B}(P_{i}^{B})$ & set of tasks which may directly follow (precede) task $i$ in backward direction, i.e., $F_{i}^{B} = \{j \in V - F_{i}^{*} - A_{i}\}$ and $P_{i}^{B} = \{j\in V|i\in F_{j}^{B}\}$ \\
        $\tau_{i j}$ & forward setup times from task $i\in V$ to task $j \in F_{i}^{F}$ \\
        $\mu_{i j}$ & backward setup times from task $i \in V$ to task $j \in F_{i}^{B}$ \\
        $\underline{\tau}_{i}$ & the smallest forward setup time from any task to task $i\in V$ \\
        $\underline{\mu}_{i}$ & the smallest backward setup time from any task to task $i\in V$ \\
        \hline
        \end{tabular}
\end{table}

Following \cite{esmaeilbeigi2016new}, SUALBP consists of
$n$ assembly tasks with precedence constraints, that require processing on $m$ different, ordered assembly stations with cycle time $c$. In the two main types of SUALBP, one of $m$ or $c$ is the objective while the maximum value of the other is fixed. 
All stations can perform all assembly tasks and the tasks assigned to a station must be completely sequenced. If a task is assigned to station $j$, all its successors must be assigned to the same or subsequent stations (i.e.,
$j, j + 1, j + 2, ... , m$). Tasks assigned to the same station
must be sequenced to satisfy the precedence constraints, if any. The deterministic processing
time of a task is provided a priori. However, the setup before a task is dependent upon the previous task in the processing sequence of the station
to which it is assigned.

For tasks assigned to a station, there are two types of sequence-dependent setups: (i) a \emph{forward setup} that captures the
setup times between any two consecutive tasks and (ii) a \emph{backward setup} that captures the setup times between the last task
of a station and the first task on the station in the next cycle. The sum of processing and the forward and backward setup times on a station constitutes the station time. The cycle time, $c$, is the maximum station time. The setups are not symmetric, i.e., the setup times from task $i$ to $j$ might be different from the setup times from task $j$ to $i$. Nevertheless, the setups satisfy the triangle inequality.

For SUALBP-1, the cycle time $c$ is fixed and the objective is to minimize the number of stations $m$. For SUALBP-2, the number of stations $m$ is fixed and the objective is to minimize the cycle time $c$. In both types, the decisions to be made are
(i) the assignment of tasks to stations; and (ii)
the sequence of the tasks assigned to each station.

We use the notation proposed by Esmaeilbeigi et al. \citeyearpar{esmaeilbeigi2016new}, as shown in the Table \ref{table:def} for both SUALBP-1 and SUALBP-2. For SUALBP-1, $\overline{c} = \underline{c} = c$. For SUALBP-2, $\overline{m} = \underline{m} = m$. To obtain all the parameters in the table, we adapt the preprocessing techniques in the literature \citep{kuroiwa2022domain, esmaeilbeigi2016new, zohali2022solving}.

\subsection{State-of-the-art MIP Model for SUALBP-1}

The state-of-the-art MIP model is the 
\emph{second station-based formulation} (SSBF) proposed by Esmaeilbeigi et al. \citeyearpar{esmaeilbeigi2016new}. Since the SSBF model can be adapted to both SUALBP-1 and SUALBP-2, we call it SSBF-1 in this section. The decision variables are:
\begin{itemize}
    \item $x_{ik}$: binary variable with value 1, iff task $i\in V$ is assigned to station $k \in FS_{i}$.
    \item $z_{i}$: integer variable for encoding the index of the station task $i \in V$ is assigned to.
    \item $u_{k}$: binary variable with value 1, iff any task is assigned to station $k$.
    \item $g_{ijk}$: binary variable with value 1, iff task $i$ is performed immediately before task $j$ on station $k$.
    \item $h_{ijk}$: binary variable with value 1, iff task $i$ is the last and task $j$ is the first task on station $k$.
    \item $r_{i}$: integer variable representing the rank of task $i$ in a sequence of all tasks. The sequence is composed of the task sequences on all the active stations, i.e., starting from the first station and ending at the last active station.
\end{itemize}

The SSBF-1 MIP model proposed by Esmaeilbeigi et al. \citeyearpar{esmaeilbeigi2016new} is as follows.
\begin{subequations}
\label{MIP1}
\begin{align}
\mathop{\min} \ & \sum_{k \in KP} u_{k} + \underline{m}  \label{6a} \qquad \qquad \quad \ \ 
    \\ 
    \text{s.t.} \ \ & \sum_{k \in FS_{i}} x_{ik} = 1, \quad
    \forall i \in V,  \label{6b}  
    \\
    & \sum_{k \in FS_{i}} k \cdot x_{ik} = z_{i}, \quad
    \forall i \in V, \label{6c} 
    \\
    & \sum_{i \in FT_{k} \cap F_{i}^{F}} g_{ijk} + \sum_{i \in FT_{k} \cap F_{i}^{B}} h_{ijk} = x_{ik},  \quad \forall i\in V, \forall k \in FS_{i},  \label{6d} 
    \\
    & \sum_{i \in FT_{k} \cap P_{j}^{F}} g_{ijk} + \sum_{i \in FT_{k} \cap P_{j}^{B}} h_{ijk} = x_{jk},  \quad \forall j\in V, \forall k \in FS_{j},  \label{6e} 
    \\
    & \sum_{i\in FT_{k}} \sum_{j\in (FT_{k} \cap F_{i}^{B})} h_{ijk} = 1, \quad \forall k \in KD,  \label{6f} 
    \\
    & \sum_{i\in FT_{k}} \sum_{j\in (FT_{k} \cap F_{i}^{B})} h_{ijk} = u_{k}, \quad \forall k \in KP,  \label{6g} 
    \\
    & r_{i} + 1 + (n - |F_{i}^{*}| - |P_{j}^{*}|) \cdot (\sum_{k \in (FS_{i} \cap FS_{j})} g_{ijk} - 1) \leq r_{j}, \quad \forall i \in V, \forall j \in F_{i}^{F}, \label{6h} 
    \\
    & r_{i} + 1 \leq r_{j}, \quad \forall (i,j) \in \mathcal{E}, \label{6i}
    \\
    & z_{i} \leq z_{j}, \quad \forall (i,j) \in \mathcal{E}, \label{6j}
    \\
    & \sum_{i \in FT_{k}} t_{i} \cdot x_{ik} + \sum_{i \in FT_{k}} \sum_{j \in (FT_{k} \cap F_{i}^{F})} \tau_{ij} \cdot g_{ijk} + \sum_{i \in FT_{k} \cap P_{i}^{B}} \mu_{ij} \cdot h_{ijk} \leq c, \quad \forall k \in KD, \label{6k}
    \\
    & \sum_{i \in FT_{k}} t_{i} \cdot x_{ik} + \sum_{i \in FT_{k}} \sum_{j \in (FT_{k} \cap F_{i}^{F})} \tau_{ij} \cdot g_{ijk} + \sum_{i \in FT_{k} \cap P_{i}^{B}} \mu_{ij} \cdot h_{ijk} \leq c \cdot u_{k}, \quad \forall k \in KP, \label{6l}
    \\
    & \sum_{i\in FT_{k} \backslash \{j\}} x_{ik} \leq (n - \underline{m} + 1) \cdot (1 - h_{jjk}), \quad \forall k \in K, \forall j \in FT_{k}, \label{6m}
    \\
    & u_{k+1} \leq u_{k}, \quad \forall k \in KP \backslash \{\overline{m}\}. \label{6s}
    \\
    & g_{ijk} \in \{0,1\}, \quad \forall k\in K, \forall i \in FT_{k}, \forall j\in (FT_{k} \cap F_{i}^{F}), \label{6n}
    \\
    & h_{ijk} \in \{0,1\}, \quad \forall k\in K, \forall i \in FT_{k}, \forall j\in (FT_{k} \cap F_{i}^{B}), \label{6o}
    \\
    & |P_{i}^{*}| + 1 \leq r_{i} \leq n - |F_{i}^{*}|, \quad \forall i \in V, \label{6p}
    \\
    & x_{ik} \in \{0,1\}, \quad \forall i \in V, \forall k \in FS_{i}, \label{6q} 
    \\
    & r_{i}, z_{i} \in \mathbb{Z}^{+}, \quad \forall i \in V, \label{6r}
\end{align}
\end{subequations}

In SSBF-1, the objective function (\ref{6a}) minimizes the number of stations. Constraint (\ref{6b}) assures that a task is assigned to a station. Constraint (\ref{6c}) connects $x_{ik}$ and $z_{i}$. Constraints (\ref{6d}) and (\ref{6e}) ensure that a task on station $k$ is followed and preceded by exactly one other task in the cyclic sequence of this station. According to constraints (\ref{6f})
and (\ref{6g}), in each cycle exactly one of the relations is a backward setup. Constraints (\ref{6h}) and (\ref{6i}) establish the precedence relations among the tasks within each
station. Note that the constraint (\ref{6h}) is inactive if tasks $i$ and $j$ are assigned to different stations. So, we add the constraint (\ref{6j}) to make sure that the precedence relations among the tasks of different stations
are satisfied. Knapsack constraints (\ref{6k}) and (\ref{6l}) guarantee that no station time exceeds
the cycle time. Constraint (\ref{6m}) ensures that only task $j$ is allocated to station $k$ when
$h_{jjk} = 1$. Constraint (\ref{6s}) assures that stations are used in the correct order and no empty station is in the middle of used stations. Constraints (\ref{6n}) to (\ref{6r}) specify domain of the variables.

Note that the decision variables $r_{i}$ and $z_{i}$ are set to continuous in \cite{esmaeilbeigi2016new}. However, doing so results in infeasible solutions being labeled as feasible for some problem instances. In addition to the MIP model, Esmaeilbeigi et al. \citeyearpar{esmaeilbeigi2016new} developed pre-processing techniques to reduce the number of variables and constraints. We implement all these techniques but omit the details here.

\subsection{State-of-the-art MIP Model for SUALBP-2}

The decision variables of the state-of-the-art SSBF-2 MIP model proposed by Zohali et al. 
\citeyearpar{zohali2022solving} are:
\begin{itemize}
    \item $x_{ik}$, $z_{i}$, $g_{ijk}$, $h_{ijk}$, $r_{i}$ as defined in Section 2.2. 
    \item $c$: continuous variable to represent the cycle time.
\end{itemize}

The SSBF-2 MIP model proposed by Zohali et al. 
\citeyearpar{zohali2022solving} is as follows.
\begin{subequations}
\label{MIP2}
\begin{align}
\mathop{\min} \ & c  \label{3a} \qquad \qquad \quad \ \ 
    \\ 
    \text{s.t.} \ \ & (\ref{6b}) - (\ref{6e}), (\ref{6h}) - (\ref{6j}), (\ref{6n}) - (\ref{6r}),
    \\
    & \sum_{i\in FT_{k}} \sum_{j\in (FT_{k} \cap F_{i}^{B})} h_{ijk} = 1, \quad \forall k \in K,  \label{3f} 
    \\
    & \sum_{i \in FT_{k}} t_{i} \cdot x_{ik} + \sum_{i \in FT_{k}} \sum_{j \in (FT_{k} \cap F_{i}^{F})} \tau_{ij} \cdot g_{ijk} + \sum_{i \in FT_{k} \cap P_{i}^{B}} \mu_{ij} \cdot h_{ijk} \leq c, \quad \forall k \in K, \label{3k}
    \\
    & \sum_{i\in FT_{k} \backslash \{j\}} x_{ik} \leq (n - \underline{m} + 1) \cdot (1 - h_{jjk}), \quad \forall k \in K, \forall j \in FT_{k}, \label{3m}
    \\
    & \sum_{i\in FT_{k}} x_{ik} \geq 1, \quad \forall k \in K, \label{3cut1}
    \\
    & c + \overline{c} \cdot (\sum_{k \in FS_{j}} k \cdot x_{jk} - \sum_{k\in FS_{i}} k \cdot x_{ik}) \geq D_{ij}, \quad \forall i \in V, j \in (F_{i}^{F} \backslash A_{i}), \label{3cut2}
    \\
    & \underline{c} \leq c \leq \overline{c}. \label{3s}
\end{align}
\end{subequations}

In SSBF-2, the objective function (\ref{3a}) minimizes the cycle time. According to constraint (\ref{3f}), in each cycle, exactly one of the relations is a backward setup. Knapsack constraint (\ref{3k}) guarantees that no station time exceeds
the cycle time. Constraint (\ref{3m}) ensures that only task $j$ is allocated to station $k$ when
$h_{jjk} = 1$.

Constraint (\ref{3cut1}) is a set of valid inequalities, based on the fact that there is at least one optimal solution with at least one task at each station because $n > m$ \citep{zohali2022solving}. Constraint (\ref{3cut2}) is another set of valid inequalities where $D_{ij}$ is a lower bound of the cycle time $c$ if tasks $i$ and $j$ are assigned to the same station \citep{esmaeilbeigi2016new}.

As above, we implement all of Zohali et al.'s pre-processing techniques but omit the details here.

\subsection{State-of-the-art Exact Approaches for SUALBP-1 and SUALBP-2}

The MIP model in Section 2.2 is the state-of-the-art exact approach for SUALBP-1 \citep{esmaeilbeigi2016new}.

For SUALBP-2, the state-of-the-art exact approach is \emph{full-featured logic-based Benders decomposition} (FFLBBD) \citep{zohali2022solving}. FFLBBD decomposes the SUALBP-2 into one master problem and a set of subproblems, one for each station. In the master problem, tasks are assigned to stations only considering the processing time and the precedence relations. In the subproblem for a station, the assigned tasks are sequenced with sequence-dependent setup times. Optimality cuts are added to the master problem based on the optimal station schedules. A relaxation of the subproblems is also included in the master problem to reduce the gap between the master problem and the subproblems. In addition, the LBBD algorithm incorporates lower and upper
bounds, pre-processing techniques, relaxations, and valid inequalities. The master problem and the subproblems are all formulated as MIP models. Given the complexity of the formulation, we refer the reader to the original paper \citep{zohali2022solving} for full details.

\subsection{Constraint Programming}

\emph{Constraint programming} (CP) is an approach to combinatorial optimization originating in the artificial intelligence community \citep{hooker2018constraint}. CP features rich variable types (e.g., interval variables and graph variables \citep{dooms2005cp}) and global constraints that represent common combinatorial sub-structure and that have an associated inference algorithm \citep{rossi2006handbook}. CP allows greater flexibility and extensibility in the types of constraints \citep{laborie2018ibm} compared to MIP, resulting in the development of many global constraints that capture modeling and inference techniques. 
However, as CP relies on local inference algorithms within global constraints to limit search and propagation of information between constraints via variable domains, a global problem view is often lacking, frequently resulting in a relatively loose CP dual bound of the global objective \citep{hooker2002logic, bockmayr1998branch, achterberg2009scip}. Nonetheless, Kuroiwa and Beck \citeyearpar{kuroiwa2023solving} showed that CP is able outperform MIP on finding and proving optimal solutions for five out of nine problem classes tested and is often able to find better quality solutions than MIP across different time limits.

CP has been used to solve SALBP-1
and SALBP-2 and outperforms MIP for larger instances \citep{bukchin2018constraint}. 
More recently, CP approaches have been applied to solving two-sided disassembly line balancing problem with AND/OR precedence 
and sequence-dependent setup times \citep{ccil2022two, kizilay2022novel}, where each station is composed of two independent workstations (sides) and precedence constraints must be addressed taking into account the predecessors of a task that might be on the other side of the assembly line, and to multi-manned SUALBP-1 where the objective is to minimize
the total number of workers required to complete pre-assigned workstation tasks \citep{guner2023constraint}. CP outperforms MIP for these non-standard SUALBP variants. In this paper, we propose CP models inspired by {\c{C}}il et al. \citeyearpar{ccil2022two}.

\subsection{Domain-Independent Dynamic Programming}

\emph{Domain-Independent Dynamic Programming} (DIDP) is a recent exact framework to solve combinatorial optimization problems \citep{kuroiwa2022domain}. The DIDP framework has been used to solve a number of problems including \emph{traveling salesman problem with time windows}, \emph{multi-commodity pickup and delivery traveling salesman problem}, and SABLP-1 \citep{kuroiwa2023solving}. It has been shown to outperform both MIP and CP models, solved by commercial MIP and CP solvers, respectively, on six of the nine problem classes tested. The success of DIDP on such problems motivates us to investigate using it to solve SUALBP.

DIDP is a model-and-solve framework for dynamic programming. A problem is represented as a dynamic programming model expressed in a domain-independent modeling language, \emph{Dynamic Programming Description Language} (DyPDL). As in MIP and CP, the model is then solved by a generic solver.

A DyPDL model is a tuple $\left \langle \mathcal{V}, \mathcal{S}^{0}, \mathcal{K}, \mathcal{T}, \mathcal{B}, \mathcal{C}, h  \right \rangle$, where
$\mathcal{V} = \{v_{1}, ..., v_{n}\}$ is the set of state variables, $\mathcal{S}^{0}$ is the target state, $\mathcal{K}$ is the set of constants, $\mathcal{T}$ is the set of transitions, $\mathcal{B}$
is the set of base cases, $\mathcal{C}$ is the set of state constraints, and $h$ is the dual bound. A solution to a DIDP model can be found by solving the recursion to determine the optimal cost of $S^{0}$. The recursion, itself, is a sequence of transitions from $S^{0}$ to a state satisfying one of the base cases.

Existing solvers for DIDP are based on heuristic state-space search approaches that have been developed over the past 50 years in artificial intelligence \citep{hart1968formal}. In particular, Kuroiwa and Beck \citeyearpar{kuroiwa2023solving} showed that \emph{complete anytime beam search} (CABS) \citep{zhang1998complete} achieves the best performance compared to other state-space search variants on each of the nine combinatorial optimization problems that were used for evaluation. Further, CABS outperformed both MIP and CP solvers (solving MIP and CP models, respectively) on six of the nine problem classes. As the name suggests, CABS, is a complete search approach: given enough time will find an optimal solution and prove optimality \citep{zhang1998complete}.

DIDP has been applied to SALBP-1 \citep{kuroiwa2022domain}, where the DIDP formulation is inspired by the state-of-the-art branch-bound-and-remember approach \citep{morrison2014application}. Our novel DIDP formulations for SUALBP are inspired by this DIDP model.

\section{CP Models for SUALBP-1 and SUALBP-2}

In this section, we present the novel CP models for SUALBP-1 and SUALBP-2.

\subsection{CP model for SUALBP-1}

The decision variables of the proposed CP model are:
\begin{itemize}
    \item $task_{i}$: interval variable \citep{laborie2018ibm} of task $i$, defined between 0 and $c$.
    \item $\overline{task}_{i}^{k}$: optional interval variable of task $i$ on station $k \in FS_{i}$, which is defined between 0 and $c$ with a size of $t_{i}$.
    \item $u_{k}$: binary variable with value 1, iff any task is assigned to station $k$.
    \item $d_{s}^{k}$: interval variable of the dummy first task on station $k$, with a size of 0.
    \item $d_{t}^{k}$: interval variable of the dummy last task on station $k$, with a size of 0.
\end{itemize}

\begin{itemize}
    \item $d_{s}^{k'}$: the next task after the dummy first task on station $k$.
    \item $d_{t}^{k'}$: the previous task before the dummy last task on station $k$.
\end{itemize}

The CP model is as follows:
\begin{subequations}
\label{CP1}
\begin{align}
\mathop{\min} \ & \sum_{k \in KP} u_{k} + \underline{m}  \label{2a} \qquad \qquad \quad \ \ 
    \\ 
    \text{s.t.} \ \  & \texttt{Alternative}(task_{i}, \overline{task}_{i}^{E_{i}}, ..., \overline{task}_{i}^{L_{i}}), \quad
    \forall i \in V,  \label{2b}  
    \\
    & \texttt{NoOverlap}(\{\overline{task}_{i}^{k} \cup \{d_{s}^{k}, d_{t}^{k}\}, \forall i \in FT_{k}\}, \{\tau_{ij}, \forall i,j \in FT_{k}\}), \quad
    \forall k \in KD \cup KP, \label{2c} 
    \\
    & \texttt{First}(\{\overline{task}_{i}^{k} \cup \{d_{s}^{k}, d_{t}^{k}\}, \forall i \in FT_{k}\}, d_{s}^{k}),  \quad \forall k \in KD \cup KP,  \label{2d} 
    \\
    & \texttt{Last}(\{\overline{task}_{i}^{k} \cup \{d_{s}^{k}, d_{t}^{k}\}, \forall i \in FT_{k}\}, d_{t}^{k}),  \quad \forall k \in KD \cup KP,  \label{2e} 
    \\
    & \texttt{PresenceOf}(\overline{task}_{i}^{k}) \leq 1, \quad \forall k \in KD, \forall i \in FT_{k},  \label{2g} 
    \\
    & \texttt{PresenceOf}(\overline{task}_{i}^{k}) \leq u_{k}, \quad \forall k \in KP, \forall i \in FT_{k},  \label{2g'} 
    \\
    & \texttt{TypeOfNext}(d_{s}^{k}) = d_{s}^{k'}, \quad \forall k \in KD \cup KP, \label{2dummy}
    \\
    & \texttt{TypeOfPrev}(d_{t}^{k}) = d_{t}^{k'}, \quad \forall k \in KD \cup KP, \label{2dummy'}
    \\
    & \texttt{StartOf}(d_{t}^{k}) + \texttt{Element}(\mu_{d_{t}^{k'}d_{s}^{k'}}) \leq c, \quad \forall k \in KD, \label{2j}
    \\
    & \texttt{StartOf}(d_{t}^{k}) + \texttt{Element}(\mu_{d_{t}^{k'}d_{s}^{k'}}) \leq c \cdot u_{k}, \quad \forall k \in KP, \label{2j'}
    \\
    & \texttt{EndBeforeStart}(\overline{task}_{i}^{k}, \overline{task}_{j}^{k}), \quad \forall k \in KD \cup DP, \forall i \in FT_{k}, \forall j \in FT_{k} \cap F_{i}^{*}, \label{2k}
    \\
    & \sum_{k \in FS_{i}} k \cdot \texttt{PresenceOf}(\overline{task}_{i}^{k}) \leq \sum_{k \in FS_{j}} k \cdot \texttt{PresenceOf}(\overline{task}_{j}^{k}), \quad \forall i \in V, \forall j \in F_{i}^{*}, \label{2l} 
    \\
    & 0 \leq task_{i} \leq c, \quad \forall i \in V, \label{2m}
    \\
    & 0 \leq \overline{task}_{i}^{k} \leq c, \quad \forall i \in V, \forall k \in FS_{i}, \label{2n}
    \\
    & 0 \leq d_{s}^{k}, d_{t}^{k} \leq c, \quad \forall k \in FS_{i}, \label{2n'}
    \\
    & u_{k} = 1, \quad \forall k \in KD, \label{2o}
    \\
    & u_{k} \in \{0,1\}, \quad \forall k \in KP. \label{2p}
\end{align}
\end{subequations}

The objective (\ref{2a}) minimizes the number of used stations. Constraint (\ref{2b}) synchronizes the fixed and optional interval variables. Constraint (\ref{2c}) ensures that the interval variables on any station form a sequence. Constraints (\ref{2d}) and (\ref{2e}) force $d_{s}^{k}$ and $d_{t}^{k}$ to be scheduled as the first and last tasks on station $k$. Constraints (\ref{2g}) and (\ref{2g'}) ensure that a task can only be assigned to an open station. Constraints (\ref{2dummy}) and (\ref{2dummy'}) guarantee that $d_{s}^{k'}$ is the next task of $d_{s}^{k}$ and $d_{t}^{k'}$ is the previous task of $d_{t}^{k}$. Constraints (\ref{2j}) and (\ref{2j'}) take the backward setup time into consideration with the help of an $\texttt{element}$ constraint. Constraint (\ref{2k}) assures that if a task and one of its followers are assigned to the same station, the follower starts after the task. Constraint (\ref{2l}) assures that a task is not assigned to a station later than the station its followers are assigned to. Constraints (\ref{2m}) to (\ref{2p}) specify variable domains.

\subsection{CP model for SUALBP-2}

The decision variables of the proposed CP model are:
\begin{itemize}
    \item $task_{i}$, $\overline{task}_{i}^{k}$, $d_{s}^{k}$, $d_{t}^{k}$, $d_{s}^{k'}$, and $d_{t}^{k'}$ as defined in Section 3.1.
    \item $c$: integer variable to represent the cycle time.
\end{itemize}

The CP model is as follows:
\begin{subequations}
\label{CP2}
\begin{align}
\mathop{\min} \ & c \label{4a} \qquad \qquad \quad \ \ 
    \\ 
    \text{s.t.} \ \  & (\ref{2b}) - (\ref{2e}), (\ref{2dummy}), (\ref{2dummy'}), (\ref{2k}), (\ref{2l}),
    \\
    & \texttt{PresenceOf}(\overline{task}_{i}^{k}) \leq 1, \quad \forall k \in K, \forall i \in FT_{k},  \label{4g} 
    \\
    & \texttt{StartOf}(d_{t}^{k}) + \texttt{Element}(\mu_{d_{t}^{k'}d_{s}^{k'}}) \leq c, \quad \forall k \in K, \label{4j}
    \\
    & 0 \leq task_{i} \leq \overline{c}, \quad \forall i \in V, \label{4m}
    \\
    & 0 \leq \overline{task}_{i}^{k} \leq \overline{c}, \quad \forall i \in V, \forall k \in FS_{i}, \label{4n}
    \\
    & 0 \leq d_{s}^{k}, d_{t}^{k} \leq \overline{c}, \quad \forall k \in FS_{i}, \label{4n'}
    \\
    & \underline{c} \leq c \leq \overline{c}. \label{4p}
\end{align}
\end{subequations}

The objective (\ref{4a}) minimizes the cycle time. Constraint (\ref{4g}) ensures that a task can only be assigned to an open station. Constraint (\ref{4j}) takes the backward setup time into consideration with the help of an $\texttt{element}$ constraint. Constraints (\ref{4m}) to (\ref{4p}) specify variable domains.

\section{DIDP Models for SUALBP-1 and SUALBP-2}

In this section, we present the novel DIDP models for SUALBP-1 and SUALBP-2.

\subsection{DIDP model for SUALBP-1}

A DIDP model can be defined using a state-transition system. A state in our DIDP model is defined by the set of unscheduled tasks, the current station that tasks are being assigned to, the first and previous (i.e., most recently assigned) tasks on the current station, and the remaining cycle time on the current station. Transitions consist of assigning a task as the next one on the current station, closing the current station to allow no further tasks to be assigned, and assigning a task as the first task on a station.

In our DIDP model, we use $\texttt{set}$ variables and $\texttt{element}$ variables, where a $\texttt{set}$ variable represents a group of elements such as customers, tasks, or items, and an $\texttt{element}$ variable represents an element of a set. DyPDL allows the use of $\texttt{resource}$ variables to capture dominance relation between states. In our DIDP model for SUALBP-1, we use an integer resource variable to represent the remaining time of the current station. If two states have the same variable values except for the resource variable, the state with a larger remaining time dominates and the other state can be soundly pruned. Similarly, in our DIDP model for SUALBP-2, we consider an integer resource variable representing the used time of the current station. Then the state with a smaller used-time value dominates another with all the other variable values being the same.

For the proposed DIDP model, we first present the state variables and the base cases of the model. We then present the recursive function in a dynamic programming form. Appendix \ref{appendix:A.1} provides an alternative perspective on this problem definition by defining the state transitions that implement the recursive function.

\noindent \emph{State variables.}
\begin{itemize}
    \item $U$: set variable for unscheduled tasks. In the target state (i.e., the initial state), $U = V$.
    \item $\kappa$: integer resource variable representing the index of the current station. In the target state, $\kappa = 0$. A smaller $\kappa$ is better. The state variable $\kappa$ is used only for computing a state-based dual bound.
    \item $p$: element variable for the previous task of the current station. In the target state, $p = d_{s}$ where $d_{s}$ is a dummy task with 0 processing time and 0 forward and backward setup with any other tasks.
    \item $f$: element variable for the first task of the current station. In the target state, $f = d_{s}$. A state keeps track of this task in order to handle the backward setup time when closing a station.
    \item $r$: integer resource variable for the remaining time (cycle time minus used time) of the current station. In the target state, $r = 0$. A larger $r$ is better.
\end{itemize}

\noindent \emph{Base case.} A base case is a set of conditions to terminate the recursion. The base case of the DIDP model is: $U = \emptyset \wedge f=d_{s}$. Note that $f=d_{s}$ is necessary since the current station has to be closed to correctly incorporate the backward setup time.

\noindent \emph{Recursive function.} We use $\mathcal{V}(U, \kappa, p, f, r)$ to represent the cost of a state. Let $U_{1} = \{j \in U | U \cap P_{j}^{*} = \emptyset\}$ be the set of tasks with all their predecessors scheduled. Let $U_{2} = \{j \in U_{1} | r \geq t_{j} + \tau_{p i}\}$ be the set of tasks in $U_{1}$ that can also be assigned to the current station, with their processing and forward setup times considered. Let $U_{3} = \{j \in U_{1} | r \geq t_{j} + \tau_{p i} + \mu_{i f}\}$ be the set of tasks in $U_{1}$ that can also be assigned to the current station, with their processing, forward setup, and backward setup times considered. The recursive function of the DIDP model is as follows:
\begin{subequations} 
\begin{align}
&\texttt{compute} \ \mathcal{V}(V, 0, d_{s}, d_{s}, 0) \label{7a} &\\  
& \mathcal{V}(U, \kappa, p, f, r) =  \label{7b} \begin{cases}
    0 & \text{ if $U=\emptyset, f=d_{s}$,} \text{\quad \ (i)}  \\
    1 + \min_{j \in U_{1}} \mathcal{V}(U\backslash\{j\}, \kappa+1, j, j, c-t_{j}) & \text{ if $U_{1} \neq \emptyset, f=d_{s}$,} \text{\quad (ii)} \\
    \min_{j \in U_{2}} \mathcal{V}(U\backslash\{j\}, \kappa, j, f, r - t_{j} - \tau_{p j}) & \text{ if $U_{2} \neq \emptyset, f \neq d_{s}$,} \text{\quad (iii)} \\
    \mathcal{V}(U, \kappa, d_{s}, d_{s}, 0) & \text{ if $U_{3} = \emptyset, f \neq d_{s}$,} \text{\quad (iv)} \\
    \infty & \text{ otherwise,} \text{\qquad \qquad (v)}
    \end{cases} \\
& U_{1} = \{j \in U | U \cap P_{j}^{*} = \emptyset\}, \quad U_{2} = \{j \in U_{1} | r \geq t_{j} + \tau_{p i} \}, \quad U_{3} = \{j \in U_{1} | r \geq t_{j} + \tau_{p i} + \mu_{i f}\} \notag \\
& \mathcal{V}(U, \kappa, p, f, r) \leq \mathcal{V}(U, \kappa^{'}, p, f, r), \quad \text{ if $\kappa \leq \kappa^{'}$ }, \label{7c} \\
& \mathcal{V}(U, \kappa, p, f, r) \leq \mathcal{V}(U, \kappa, p, f, r^{'}), \quad \text{ if $r \geq r^{'}$ }, \label{7d} \\
& \mathcal{V}(U, \kappa, p, f, r) \geq \max 
 \label{7e} 
\begin{cases}
    \ceil[\Big]{ \frac{\underline{\mu}_{f} +  \sum_{i \in U} (\underline{\tau}_{i} + t_{i}) - (\overline{m} - \kappa) \cdot (\max_{i \in U}{\underline{\tau}_{i}}) + \max(\underline{m} - \kappa, 0) \cdot (\min_{i \in U}{\underline{\mu}_{i}}) - r}{c} }, \text{\;  (i)} \\
    \ceil[\Big]{ \frac{\underline{\mu}_{f} +  \sum_{i \in U} t_{i} - r}{c} }, \text{\qquad \qquad \qquad \qquad \qquad \qquad \qquad \qquad \quad \ (ii)}  \\
    \ceil{ \frac{\sum_{i\in U} t_{i} - r}{c} }, \text{\qquad \qquad \qquad \qquad \qquad \qquad \qquad \qquad \quad \quad \; (iii)}  \\
    \sum_{i \in U} w_{i}^{2} + \ceil{ \sum_{i\in U} w_{i}^{'2} - l^{2} }, \text{\qquad \qquad \qquad \qquad \qquad \quad \quad (iv)} \\
    \ceil{ \sum_{i\in U} w_{i}^{3} - l^{3} }. \text{\qquad \qquad \qquad \qquad \qquad \qquad \qquad \qquad \quad \ (v)}
\end{cases}
\end{align}
\end{subequations} 

The term (\ref{7a}) is to compute the cost of the target state. Equation (\ref{7b}) is the main recursion of the DIDP model. Specifically, (5b-i) refers to the base case, while (5b-ii) corresponds to assigning the first task to the current station. The recursion here can be understood as assigning the cost of the current state as equal to one more than the state where $j$ is the first task on station $\kappa + 1$. Case (5b-iii) represents assigning the next task after other tasks have already been assigned to the current station. Since the number of used stations is not increased, the cost of the current state is the same as the cost of the state where $j$ is scheduled next on station $\kappa$. Note that task $j$ is assigned as the new last task on station $\kappa$. Cases (5b-iv) and (5b-v) correspond to closing the current station and detecting dead-ends, respectively. Inequalities (\ref{7c}) and (\ref{7d}) formulate state domination in two scenarios: if other variables are equal a state with smaller $\kappa$ or larger remaining time dominates. Term (\ref{7e}) states the state-based dual bounds. Both (5e-1) and (5e-ii) are novel dual bounds whose validity we prove below. (5e-iii), (5e-iv), and (5e-v) were first used for SALBP-1 \citep{kuroiwa2022domain} and are also valid here.

The numeric constants $w^{2}, w^{'2}, w^{3}$ are indexed by a task $i$ and depend on $t_{i}$, as shown in Table \ref{table:knapsack_dual}. These values are obtained by ignoring precedence relations and were originally proposed by Scholl and Klein \citeyearpar{scholl1997salome}.

\begin{table}[tp]
    \caption{Numeric constants for calculating a knapsack-based dual bound.}
    \label{table:knapsack_dual}
    \centering
    \begin{tabular}{c|ccc||c|ccccc}
        \toprule[1.0pt]
        $t_{i}$ & (0, c/2) & c/2 & (c/2, c] & $t_{i}$ & (0, c/3) & c/3 & (c/3, c/2) & 2c/3 & (2c/3, c] \\
        \midrule[1.0pt]
        $w_{i}^{2}$ & 0 & 0 & 1 & $w_{i}^{3}$ & 0 & 1/3 & 1/2 & 2/3 & 1 \\
        $w_{i}^{'2}$ & 0 & 1/2 & 0 & & \\
        \bottomrule[1.0pt]
    \end{tabular}
\end{table}

\subsection{Correctness of Novel Dual Bounds}

For the proposed state-based dual bound (5e-i), $\underline{\tau}_{i}$ is the smallest forward setup time from any task to task $i$ and $\underline{\mu}_{i}$ is the smallest backward setup time from any task to task $i$. The value $\max_{i \in U}{\underline{\tau}_{i}}$ is the largest minimum forward setup time to any unscheduled task. Similarly, $\min_{i \in U}{\underline{\mu}_{i}}$ is the smallest minimum backward setup time to any unscheduled task.

\noindent \begin{theorem}
\label{bound1a}
Term (5e-i) is a valid lower bound of the number of additional stations to be used at the current state.
\end{theorem}

\proof{Proof.} 

\begin{figure}[bp]
  \centering
  \includegraphics[width=0.35\textwidth]{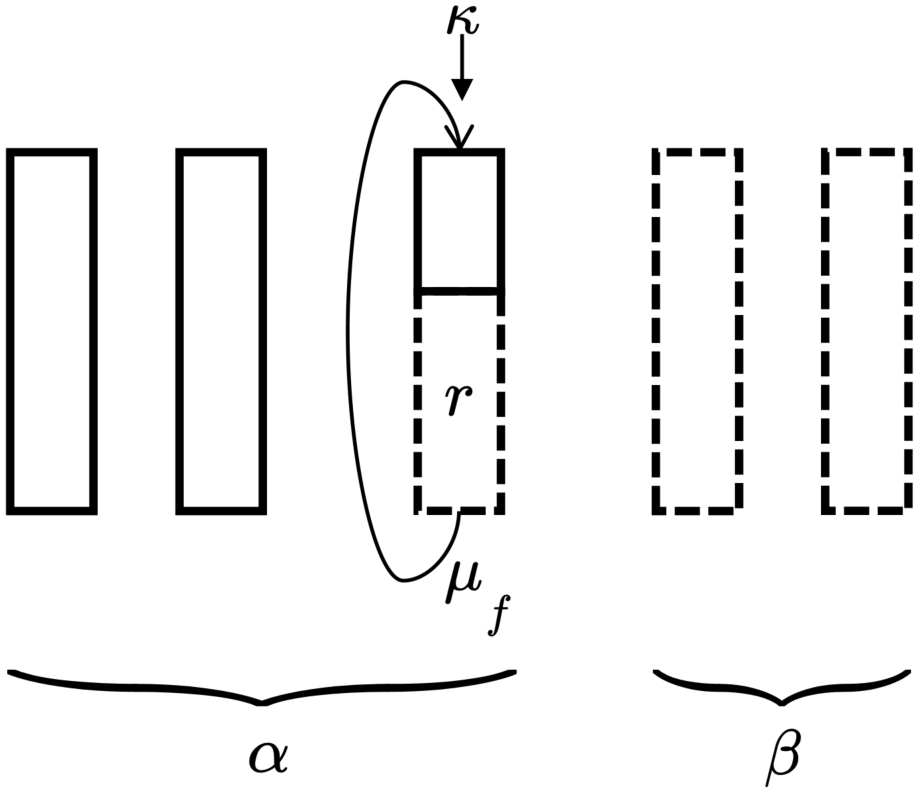}
  \caption{Illustration of dual bound (5e-i).}
  \label{proof1}
\end{figure}

In Fig. \ref{proof1} solid rectangles are stations that have tasks assigned, while dashed rectangles are stations that have no tasks assigned yet. Let the number of remaining stations to be used be $\beta$, the number of stations that are already used be $\alpha$, and the current station, $\kappa$. $r$ contributes to the total time of the current station $\kappa$ and is hence excluded from the calculation of $\beta$. Let the total time represented by the dashed rectangles in the best possible solution given the current state be $T$, then,
\begin{equation}
    T = \underline{\mu}_{f} + \sum_{i\in U}t_{i} + \sum_{i\in U}\delta_{i}^{*}
\end{equation}
where $\delta_{i}^{*}$ is either the forward or the backward setup time starting from task $i$ (i.e., $\tau_{i j}$ where $j$ is the immediate successor of $i$ in the station that $i$ is assigned to, or $\mu_{i j}$ where $j$ is the first task and $i$ is the last task in the station that $i$ is assigned to) in the best possible solution reachable from the current state. Without knowing the best possible solution reachable from the current state, the value of $\sum_{i\in U}\delta_{i}^{*}$ is unclear. A lower bound can be obtained as
\begin{equation} \label{8}
    \sum_{i\in U}\delta_{i}^{*} \geq \sum_{i\in U}\underline{\tau}_{i} - (\overline{m} - \kappa) \cdot \max_{i \in U}{\underline{\tau}_{i}} + \max(\underline{m} - \kappa, 0) \cdot \min_{i \in U}{\underline{\mu}_{i}}.
\end{equation}
The inequality (\ref{8}) is true since: (i) each station needs to consider a backward setup time and for the $\underline{m} - \kappa$ stations that are definitely used in the future, $(\underline{m} - \kappa) \cdot \min_{i \in U}{\underline{\mu}_{i}}$ is a lower bound on the real backward setup times on these stations,\footnote{Recall that $\underline{m}$ is a lower bound on the number of machines uses.} however, since it is possible that $\underline{m} < \kappa$ at some states and negative backward setup times should not be added,  $\max(\underline{m} - \kappa, 0) \cdot \min_{i \in U}{\underline{\mu}_{i}}$ is used; and (ii) as the remaining setup times are all forward, $\sum_{i\in U}\underline{\tau}_{i} - (\overline{m} - \kappa) \cdot \max_{i \in U}{\underline{\tau}_{i}}$ is a lower bound of all the forward setup times. Thus,
\begin{equation}
    T \geq \underline{\mu}_{f} + \sum_{i\in U}t_{i} + \sum_{i\in U}\underline{\tau}_{i} - (\overline{m} - \kappa) \cdot \max_{i \in U}{\underline{\tau}_{i}} + \max(\underline{m} - \kappa, 0) \cdot \min_{i \in U}{\underline{\mu}_{i}}.
\end{equation}
Given that $\beta = \ceil{\frac{T-r}{c}}$, we finally have
\begin{equation}
    \beta \geq \ceil[\Bigg]{ \frac{\underline{\mu}_{f} +  \sum_{i \in U} (\underline{\tau}_{i} + t_{i}) - (\overline{m} - \kappa) \cdot (\max_{i \in U}{\underline{\tau}_{i}}) + \max(\underline{m} - \kappa, 0) \cdot (\min_{i \in U}{\underline{\mu}_{i}}) - r}{c} }.
\end{equation}

A corner case is the target state where $\kappa = r = 0$ and $f = d_{s}$. In that case, the inequality is simplified to 
\begin{equation}
    \beta \geq \ceil[\Bigg]{ \frac{ \sum_{i \in U} (\underline{\tau}_{i} + t_{i}) + \underline{m} \cdot (\min_{i \in U}{\underline{\mu}_{i}}) - \overline{m} \cdot (\max_{i \in U}{\underline{\tau}_{i}}) }{c} },
\end{equation}
which is also valid as $\sum_{i \in U} (\underline{\tau}_{i} + t_{i}) + \underline{m} \cdot (\min_{i \in U}{\underline{\mu}_{i}}) - \overline{m} \cdot (\max_{i \in U}{\underline{\tau}_{i}})$ is smaller than the real total station time of all active stations.
\myright{\Halmos} \endproof

Similarly, if forward setup times are not considered and only the backward setup of the current station is counted, (5e-ii) is obtained. We omit the proof of its correctness. In most of cases, (5e-i) is greater than (5e-ii), but there are exceptions. We present the proof of the case where (5e-ii) is greater than or equal to (5e-i) as follows:

\begin{theorem}
\label{bound1b}
Term (5e-i) does not dominate (5e-ii).
\end{theorem}

\proof{Proof.} 
Let $\overline{m} = \kappa + 2$ and $\kappa > \underline{m}$. Also, let $U = \{i,j\}$,  $\underline{\tau}_{i} = 2$, and $\underline{\tau}_{j} = 12$. Then (5e-i) becomes less than or equal to (5e-ii) as follows:
\begin{equation*}
    \ceil[\Big]{ \frac{\underline{\mu}_{f} +  \sum_{i \in U} t_{i} + 2 + 12 - 2 \times 12 + 0 - r}{c} } = \ceil[\Big]{ \frac{\underline{\mu}_{f} +  \sum_{i \in U} t_{i} - 10 - r}{c} } \leq \ceil[\Big]{ \frac{\underline{\mu}_{f} +  \sum_{i \in U} t_{i} - r}{c} }.
\end{equation*}
Let $\underline{\mu}_{f}=1, \sum_{i \in U} t_{i} = 50, r=1$, 
and $c=10$, then at this state, (5e-i)$= 4 < 5 =$(5e-ii).
Thus, (5e-i) does not dominate (5e-ii).
\endproof

\subsection{DIDP model for SUALBP-2}

This model is similar to the DIDP model of SUALBP-1, with a new state variable $c$ to represent the cycle time and the replacement of $r$ by $t^{c}$. In this model, we use $t^{c}$ to keep track of the accumulated time on the current station. By contrast, in the model of SUALBP-1, we use $r$ to keep track the remaining time of the current station. By definition $t^{c} + r = c'$ where $c'$ is the cycle time of the current station. As above, the transition-centric view of the model definition is provided in Appendix \ref{appendix:A.2}.

\noindent \emph{State variables.}
\begin{itemize}
    \item $U, p, f$: these set variables are the same as the DIDP model for SUALBP-1.
    \item $\kappa$: integer variable for the current station. In the target state, $\kappa = 0$. Since the objective of SUALBP-2 is not to minimize the number of active stations, a smaller number of active stations is not preferred anymore, and so $\kappa$ is not a resource variable.
    \item $t^{c}$: integer resource variable for the used time of the current station. In the target state, $t^{c} = 0$ represents that no task is assigned to the current station. A smaller $t^{c}$ is better.
    \item $c$: integer resource variable for the cycle time. A smaller $c$ is better. In the target state, $c=0$. Note that the state variable $c$ is used only for computing a state-based dual bound.
\end{itemize}

\noindent \emph{Base case.} The base case of the DIDP model is: $U = \emptyset \wedge f=d_{s}$. Note that $f=d_{s}$ is necessary since in the base case one has to include the backward setup time for the last station.

\noindent \emph{Recursive function.} We use $\mathcal{V}(U, \kappa, p, f, t^{c}, c)$ to represent the cost of a state. Let $U_{1} = \{j \in U | U \cap P_{j}^{*} = \emptyset\}$. The recursive function of the DIDP model is as follows:
\begin{subequations} 
\begin{align}
&\texttt{compute} \ \mathcal{V}(V, 0, d_{s}, d_{s}, 0, 0) \label{8a} &\\  
& \mathcal{V}(U, \kappa, p, f, t^{c}, c) = \label{8b} \\
& \qquad \begin{cases}
    0 & \text{if $U=\emptyset, f=d_{s}$,}  \notag \qquad \quad \ \ \text{(i)} \\
    \min_{j \in U_{1}} ( \max(t_{j}, \mathcal{V}(U\backslash\{j\}, \kappa+1, j, j, t_{j}, \max(c, t_{j})))) & \text{if $U_{1} \neq \emptyset, f=d_{s}, \kappa < m$,} \ \text{(ii)} \\
    \min_{j \in U_{1}} (\max(t^{c} + t_{j} + \tau_{p j}, \mathcal{V}(U\backslash\{j\}, \kappa, j, f,  \notag \\ 
    \qquad \qquad \qquad t^{c} + t_{j} + \tau_{p j}, \max(c, t^{c} + t_{j} + \tau_{p j})))) & \text{if $U_{1} \neq \emptyset, f \neq d_{s}$,} \qquad \quad \ \text{(iii)} \\
    \max(t^{c} + \mu_{p f}, \mathcal{V}(U, \kappa, d_{s}, d_{s}, t^{c} + \mu_{p f}, \max(c, t^{c} + \mu_{p f}))) & \text{if $U_{1} = \emptyset, f \neq d_{s}$,} \qquad \quad \ \text{(iv)} \\
    \infty & \text{otherwise,} \qquad \qquad \qquad \ \text{(v)}
    \end{cases} \\ 
& U_{1} = \{j \in U | U \cap P_{j}^{*} = \emptyset\}, \notag \\
& \mathcal{V}(U, \kappa, p, f, t^{c}, c) \leq \mathcal{V}(U, \kappa, p, f, t^{c'}, c), \quad \text{ if $t^{c} \leq t^{c'}$ }, \label{8c} \\
& \mathcal{V}(U, \kappa, p, f, t^{c}, c) \leq \mathcal{V}(U, \kappa, p, f, t^{c}, c^{'}), \quad \text{ if $c \leq c^{'}$ }, \label{8d} \\
& \mathcal{V}(U, \kappa, p, f, t^{c}, c) \geq \max \label{8e} 
\begin{cases}
    \ceil[\Big]{ \frac{ \sum_{i \in U} (\underline{\tau}_{i} + t_{i}) + t^{c} + (m - \kappa) \cdot (\min_{i \in U}{\underline{\mu}_{i}} - \max_{i \in U}{\underline{\tau}_{i}}) + \underline{\mu}_{f}}{\min (m, m - \kappa + 1)} - c}, \qquad \text{(i)} \\
    \ceil[\Big]{ \frac{ \sum_{i \in U} t_{i} + t^{c} }{\min (m, m - \kappa + 1)} - c}. \qquad \qquad \qquad \qquad \qquad \qquad \qquad  \text{(ii)}
\end{cases}
\end{align}
\end{subequations} 

The term (\ref{8a}) is to compute the objective of the target state. Equation (\ref{8b}) is the main recursion of the DIDP model. Specifically, (11b-i) handles the base cases, while (11b-ii) refers to assigning the first task to the current station. The recursion here can be understood as assigning the cost of the current state as equal to the maximum station time of the state where $j$ is the first task on station $\kappa + 1$. Case (11b-iii) corresponds to assigning the next task to a station after other tasks have been assigned. Cases (11b-iv) and (11b-v) deal with closing the current station and dead-ends, respectively. Inequality (\ref{8c}) and (\ref{8d}) describe the state domination in two scenarios: (i) smaller used time dominates larger one and (ii) smaller cycle time dominates larger one, given all the other state variables are the same. Term (\ref{8e}) formulates two novel state-based state-based dual bounds.

\subsubsection{Correctness of Novel Dual Bounds}

\noindent \begin{theorem}
\label{bound2a}
Term (11e-i) is a valid dual bound of the cycle time at the current state.
\end{theorem}

\begin{figure}[bp]
  \centering
  \includegraphics[width=0.4\textwidth]{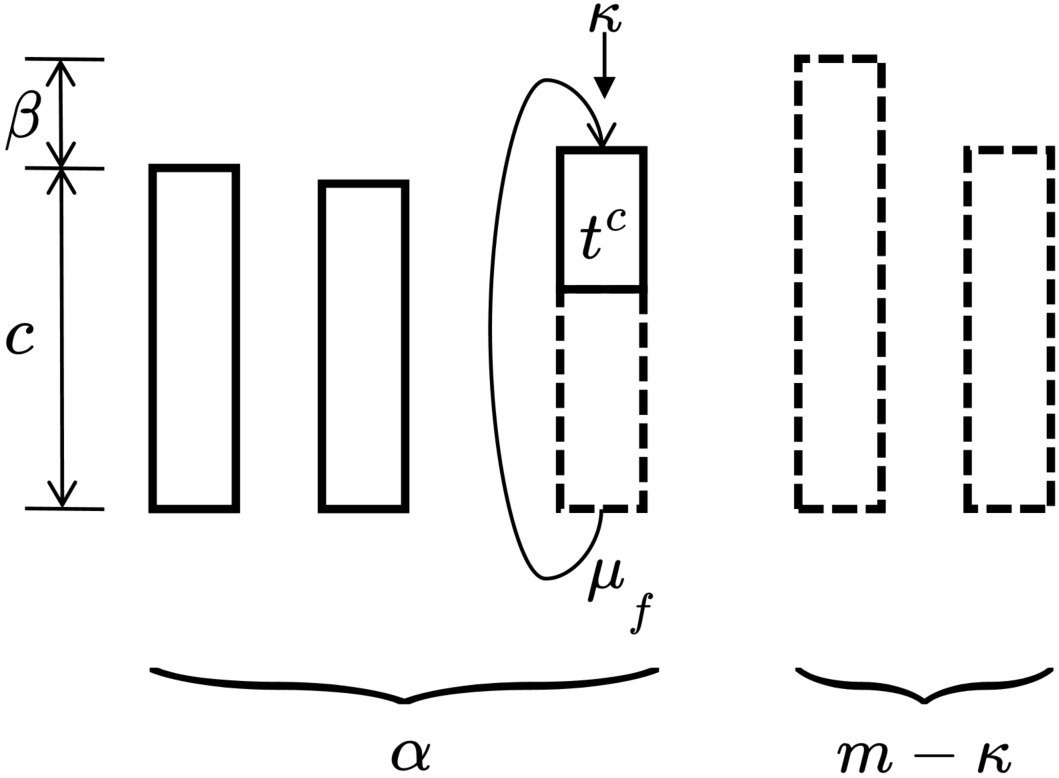}
  \caption{Illustration of dual bound (11e-i).}
  \label{proof2}
\end{figure}

\proof{Proof.} 
In Fig. \ref{proof2}, solid rectangles are stations that have tasks assigned, while dashed rectangle are stations that have no tasks assigned yet. Let the number of stations that are already used be $\alpha$, the current station be $\kappa$, the number of remaining stations to be used be $m -\kappa$, the current cycle time be $c$, and the difference between the future larger cycle time and $c$ be $\beta$. $t^{c}$ contributes to the total time of the current station $\kappa$ and is hence excluded from the calculation of $c$. Let the total time represented by the dashed rectangles in the best possible solution given the current state be $T$, then,
\begin{equation}
    T = \underline{\mu}_{f} + \sum_{i\in U}t_{i} + \sum_{i\in U}\delta_{i}^{*}.
\end{equation}
where $\delta_{i}^{*}$ is either the forward or the backward setup time starting from task $i$ as defined in the proof of Theorem 1. Similar to the proof of Theorem 1, the following inequality is valid.
\begin{equation}
    \sum_{i\in U}\delta_{i}^{*} \geq \sum_{i\in U}\underline{\tau}_{i} - (m - \kappa) \cdot \max_{i \in U}{\underline{\tau}_{i}} + (m - \kappa) \cdot \min_{i \in U}{\underline{\mu}_{i}}.
\end{equation}
Thus,
\begin{equation}
    T \geq \underline{\mu}_{f} + \sum_{i\in U}t_{i} + \sum_{i\in U}\underline{\tau}_{i} - (m - \kappa) \cdot \max_{i \in U}{\underline{\tau}_{i}} + (m - \kappa) \cdot \min_{i \in U}{\underline{\mu}_{i}}.
\end{equation}
Given that $\beta = \ceil{\frac{T+t^{c}}{m-\kappa+1} - c}$, we finally have
\begin{equation}
    \beta \geq \ceil[\Bigg]{ \frac{ \sum_{i \in U} (\underline{\tau}_{i} + t_{i}) + t^{c} + (m - \kappa) \cdot (\min_{i \in U}{\underline{\mu}_{i}} - \max_{i \in U}{\underline{\tau}_{i}}) + \underline{\mu}_{f}}{(m - \kappa + 1)} - c}.
\end{equation}
Also, $\kappa \leq m$ is always true, as the only case that can increase $\kappa$ requires $\kappa < m$, as shown in (11b-ii). Thus, $m-\kappa + 1 \geq 1$ and we do not need to handle 0 divisors in the dual bounds.

Another corner case is the target state where $\kappa = t^{c} = c = 0$ and $f = d_{s}$. In that case, the inequality is simplified to 
\begin{equation}
    \beta \geq \ceil[\Bigg]{ \frac{ \sum_{i \in U} (\underline{\tau}_{i} + t_{i}) + m \cdot (\min_{i \in U}{\underline{\mu}_{i}} - \max_{i \in U}{\underline{\tau}_{i}}) }{m}}.
\end{equation}
Therefore, we finally get 
\begin{equation}
    \beta \geq \ceil[\Bigg]{ \frac{ \sum_{i \in U} (\underline{\tau}_{i} + t_{i}) + t^{c} + (m - \kappa) \cdot (\min_{i \in U}{\underline{\mu}_{i}} - \max_{i \in U}{\underline{\tau}_{i}}) + \underline{\mu}_{f}}{\min (m, m - \kappa + 1)} - c}.
\end{equation}
\myright{\Halmos} \endproof

Similarly, (11e-ii) is obtained if forward and backward setup times are not considered and we omit the proof of its correctness. In majority of cases, (11e-i) leads to a larger dual bound than (11e-ii). Nevertheless, we prove the necessity of (11e-ii) by displaying a case where (11e-ii) needs to be strictly greater than (11e-i).

\begin{theorem}
\label{bound2}
Term (11e-i) does not dominate (11e-ii).
\end{theorem}

\proof{Proof.} 
Let $m = \kappa + 2$, $U = \{i,j\}$,  $\underline{\tau}_{i} = 2$, and $\underline{\tau}_{j} = 12$. Also, let $\underline{\mu}_{i} = 1$, $\underline{\mu}_{j} = 1$, and $\underline{\mu}_{f} = 1$. Then (11e-i) becomes less than or equal to (11e-ii) as follows:
\begin{equation*}
    \ceil[\Big]{ \frac{\sum_{i \in U} t_{i} + 2 + 12 + 2 \times (1 - 12) + 1 + t^{c}}{\min (m, 2 + 1)} - c} = \ceil[\Big]{ \frac{\sum_{i \in U} t_{i} - 7 + t^{c}}{\min (m, 3)} - c} \leq \ceil[\Big]{ \frac{ \sum_{i \in U} t_{i} + t^{c} }{\min (m, 3)} - c}.
\end{equation*}
Let $\sum_{i \in U} t_{i} = 30, t^{c}=1, m=5$, and $c=5$, then at this state, (11e-i)$=3<6=$(11e-ii).
Thus, (11e-i) does not dominate (11e-ii).
\endproof

\section{Local Improvement for SUALBP-2}

In this section, we present a local improvement algorithm based on the DIDP model for SUALBP-2. The idea is simple: given a feasible solution of SUALBP-2, the task schedule of each station can be independently re-optimized to reduce the station time, which might lead to a better cycle time. We call this process the ``local improvement".

The local improvement method can be done combined with any anytime search algorithm. If the original algorithm is exact, it remains exact after the combination. Whenever a new incumbent with objective $c$ is found for SUALBP-2, we try to refine this solution. Specifically, we first sort the stations in a decreasing order of station time. Then starting from the first station, we search for the sequence that minimizes station time. If no improvement exists on a station, the local improvement algorithm cannot (further) decrease the cycle time and it exits. If a better solution with objective $c'$ is hence found, we return to where the SUALBP-2 algorithm paused with $c'$ as the incumbent cost.

The \emph{local improvement} (LI) algorithm is shown in Algorithm \ref{algo:LI}. The initialization of the algorithm is conducted in line 1, where the stations are sorted in a decreasing order of the station times in the incumbent solution $\mathbf{X}$. $K$ is the set of sorted stations while $c_{k}$ is the station time of station $k \in K$.
In the loop from line 2 to line 5, each station in $K$ is popped in order and the corresponding station time is re-optimized. If the station time is not reduced or the reduced station time is greater than or equal to the station time of the next station in $K$, the local improvement is finished. The entire solution of SUALBP-2 is then updated and returned as shown in line 6 and line 7.

\begin{algorithm}[bp]
    \caption{Local Improvement}
    \label{algo:LI}
    \BlankLine
    \KwIn{$\mathbf{X}$ - incumbent solution of the original problem}
    
    $K, \{c_{k}, \forall k\in K\} \leftarrow \texttt{StationTimeDecreasing}(\mathbf{X})$;

    \Repeat{$c_{k}^{'} \geq c_{k}$ \emph{or} $c_{k}^{'} \geq c_{k+1}$}{
        $k \leftarrow \texttt{pop}(K)$;
        
        $\mathbf{X}_{k}^{'}, c_{k}^{'} \leftarrow \texttt{LocalImprovement}(\mathbf{X}_{k})$;
    }          

    $\mathbf{X}^{'} \leftarrow \texttt{UpdateSolution}(\mathbf{X, \mathbf{X}_{k}^{'}})$;
    
    \textbf{return} $\mathbf{X}^{'}, c_{k}^{'}$;
\end{algorithm}

The task re-sequencing subproblem on a station is very similar to the travelling salesman problem with precedence constraints, which can be solved with DIDP \citep{kuroiwa2023solving}. Thus, we formulate each of the re-sequencing problem as a DIDP model. Since the DIDP model is simple, we present it in Appendix \ref{appendix:A.3}.

\section{Experiments}

In this section, we compare the performance of
our DIDP and CP models against the state-of-the-art MIP and FFLBBD models \citep{esmaeilbeigi2016new, zohali2022solving} on the SBF2 data set (788 instances) \citep{scholl2013assembly}.\footnote{\url{https://assembly-line-balancing.de/sualbsp/data-set-of-scholl-et-al-2013/}}

The data set considers four
levels of setup times, determined by parameter $\alpha$ that specifies the ratio of the average setup time to the average task processing time. Higher values of $\alpha$ represent larger setup
times. Zohali et al. \citeyearpar{zohali2022solving} create four different data
sets by setting $\alpha$ values equal to 0.25, 0.50, 0.75, and
1.00.

The SBF2 data set was designed for SUALBP-1 and provides the optimal number of stations for type-1. However, it does not provide the number of stations nor the optimal cycle time for type-2. Zohali et al. \citeyearpar{zohali2022solving} consider
$m = \sum_{i \in V} t_{i}/c$ where $c$ is the pre-defined cycle time obtained from the SBF2 data set for each instance. Thus, this data set can also be used for testing SUALBP-2.

Zohali et al. \citeyearpar{zohali2022solving} cluster these instances into four classes:

\begin{itemize}
    \item Data set A: small (132 instances) with up to 25 tasks.
    \item Data set B: medium (140 instances) with 28 to 35 tasks.
    \item Data set C: large (188 instances) with 45 to 70 tasks.
    \item Data set D: extra-large (328 instances) with 75 to
111 tasks.
\end{itemize}

For evaluation, we use the fraction of instances that are solved and proved optimal over time and the fraction of instances over primal integral \citep{berthold2013measuring}. The primal integral considers the balance between the solution quality and computational time. For an optimization problem, let $s^{t}$ be a solution found at time $t$ by an algorithm, $s^{*}$ be an optimal (or best-known) solution, and $c$ be a function mapping a solution to its corresponding cost. The primal gap function $p$ is defined as follows:
\begin{equation}
    p(t) =   
    \begin{cases}
        1   & \text{if no } s^{t} \text{ or } c(s^{t})c(s^{*}) < 0, \\
        0   & \text{if } c(s^{t}) = c(s^{*}) = 0, \\
        \frac{|c(s^{*}) - c(s^{t})|}{max\{|c(s^{*})|, |c(s^{t})|\}}    & \text{otherwise}.
    \end{cases}
\end{equation}

The primal gap takes a value in range $[0, 1]$, with lower values being better. We use $p(T)$, the
primal gap at the time limit $T$, as a metric of the final solution
quality. Let $t_{i} \in [0, T]$ for $i = 1, ..., L - 1$ be a time point when a new incumbent solution is found by an algorithm, $t_{0} = 0$, and $t_{L} = T$. The primal integral is defined as:
\begin{equation}
    P(T) = \sum_{i=1}^{L}p(t_{i}-1) \cdot (t_{i} - t_{i-1}).
\end{equation}

The primal integral takes a value in $[0, T]$, with lower values being better. $P(T)$ decreases
if a better solution is found with the same computational time or the same solution cost is achieved faster. When an
instance is proved to be infeasible at time $t$, we use $p(t) = 0$,
corresponding to the time to prove infeasibility.

For the DIDP models, we use the state-of-the-art solver based on complete anytime beam search (CABS) \citep{kuroiwa2023solving} in didp-rs v0.4.0.\footnote{https://didp.ai/} We have also tested the local improvement algorithm combined with DIDP models for the original problem and subproblems of SUALBP-2. Since we use CABS as the DIDP solver, we call the entire algorithm \emph{local improvement CABS} (LICABS). 
For the CP models, we use CP Optimizer 20.1.0 \citep{cplex}. For the MIP models, we use Gurobi 9.5.1 \citep{gurobi}. All the experiments are implemented in Python 3.8. Each instance is run for 1800 seconds on a single thread on a Ubuntu 22.04.2 
LTS machine with Intel Core i7 CPU and 16 GB memory.

\subsection{Results for SUALBP-1}

\begin{figure}[tp]
\centering
\begin{minipage}[t]{0.48\textwidth}
\centering
  \includegraphics[width=1.00\textwidth]{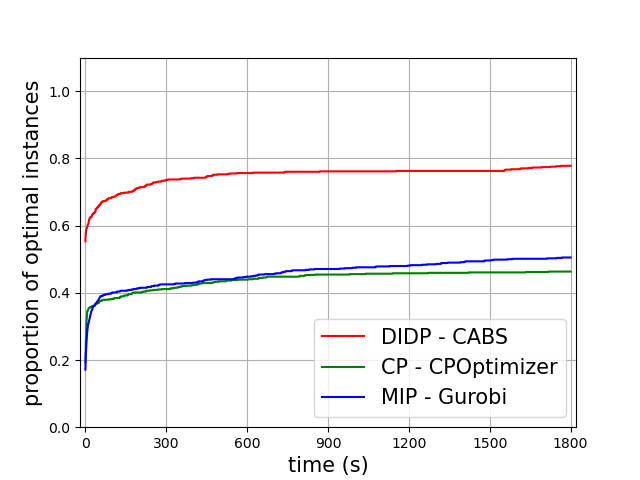}
  \caption{Ratio of instances solved to optimality over time for SUALBP-1}
  \label{time1}
\end{minipage}
\hspace{0.1cm}
\begin{minipage}[t]{0.48\textwidth}
\centering
  \includegraphics[width=01.00\textwidth]{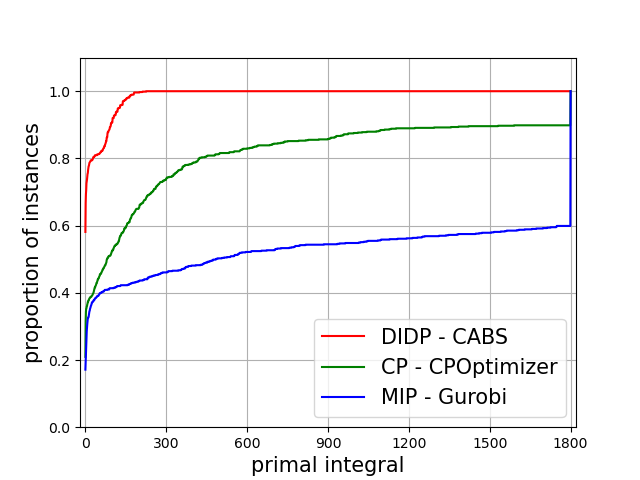}
  \caption{Ratio of instances over primal integral for SUALBP-1}
  \label{pi1}
\end{minipage}
\end{figure}

The results on SUALBP-1 are shown in Fig. \ref{time1} and \ref{pi1}. Better performance is indicated by curves closer to the top left corner of the graphs. The DIDP model outperforms CP and MIP models. More specifically, Fig. \ref{time1} shows that the DIDP model finds optimal solutions and proves optimality for more instances in a shorter computation time than the CP and MIP models. In 1 second, DIDP finds and proves optimality on 55\% of the instances. MIP and CP cannot achieve the same performance in 1800 seconds. At the 1800 second time limit, DIDP has found and proved optimality for 78\% of the problem instances compared to 45\% and 50\% for CP and MIP, respectively.

Fig. \ref{pi1} shows that DIDP also finds high-quality solutions faster than CP and MIP models. As is often observed, although the CP model proves fewer optimal instances than the MIP model, it demonstrates a higher solution quality \citep{hooker2018constraint}. However, DIDP outperforms CP and MIP on both measures.

\subsection{Results for SUALBP-2}

\begin{figure}[tp]
\centering
\begin{minipage}[t]{0.48\textwidth}
\centering
  \includegraphics[width=1.00\textwidth]{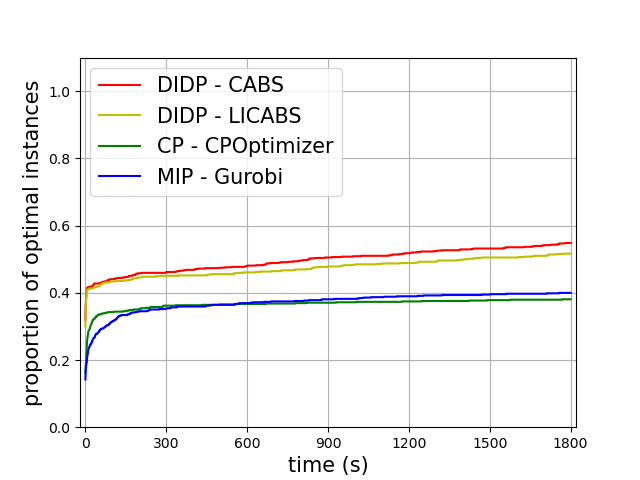}
  \caption{Ratio of instances solved to optimality over time for SUALBP-2}
  \label{time2}
\end{minipage}
\hspace{0.1cm}
\begin{minipage}[t]{0.48\textwidth}
\centering
  \includegraphics[width=01.00\textwidth]{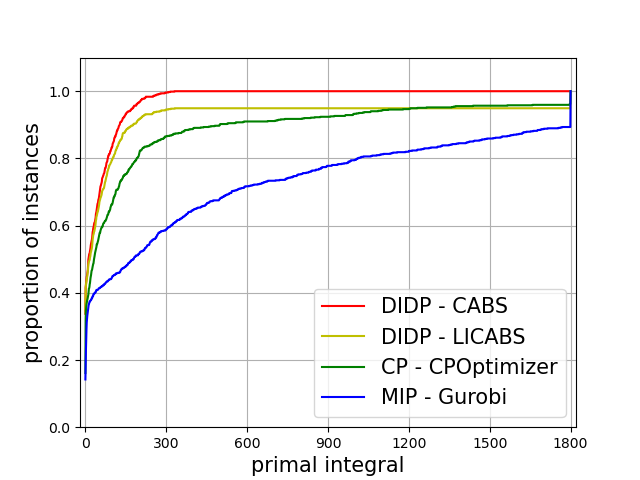}
  \caption{Ratio of instances over primal integral for SUALBP-2}
  \label{pi2}
\end{minipage}
\end{figure}

The results on SUALBP-2 are shown in Fig. \ref{time2} and \ref{pi2}. The LICABS results are discussed in Section 6.4. As above, better performance is indicated by curves closer to the top and left-side of the graphs. In general, the two graphs demonstrate similar though reduced performance of the three models compared to SUALBP-1. The DIDP model continues to substantially outperform CP and MIP models. Fig. \ref{time2} shows that the DIDP model finds optimal solutions and proves optimality for more instances in a shorter computation time than CP and MIP models. 
In 1 second, DIDP finds and proves optimality on 32\% of the instances. CP achieves the same level at 30 seconds and is by that measure 30 times slower. MIP reaches the same level at 120 seconds. At the 1800 second time limit, DIDP has found and proved optimality for 55\% of the problem instances compared to 37\% and 41\% for CP and MIP, respectively.

Fig. \ref{pi2} shows that DIDP finds high-quality solutions faster than CP and MIP models. Although the CP model proves fewer optimal instances than the MIP model, it demonstrates a higher solution quality than the MIP model and is almost equal to the DIDP model at the time limit. However, DIDP achieved the same level of quality as CP and MIP in considerably less than 180 and 140 seconds, respectively.

Overall the experimental results for SUALBP-1 and SUALBP-2 show that DIDP substantially outperforms MIP and CP while CP is competitive with MIP: CP finds better solutions than MIP in a given time but trails MIP in the number of solutions proved optimal over time.

\subsection{Comparison to Zohali et al. for SUALBP-2}

As noted, the state-of-the-art exact algorithm for SUALBP-2 is the \emph{full-featured logic based Benders decomposition} (FFLBBD) approach proposed by Zohali et al. \citeyearpar{zohali2022solving}. We were not able to obtain their source code but, as they ran experiments on the same SBF2 data set, we can compare the quality of results to what the authors provided in the online appendix to that paper. The detailed comparison is shown in Table \ref{sota_no_better}. The results of FFLBBD are taken from paper by Zohali et al. \citeyearpar{zohali2022solving} with different hardware, but the same time and memory limits of 1800 seconds and 16 GB. For the overall comparison, DIDP proves 438 instances to optimality while FFLBBD only proves 371. 
Since these benchmark instances are first investigated by Zohali et al. \citeyearpar{zohali2022solving}, according to their results, FFLBBD reports solving 9 instances to optimality that DIDP does not, while DIDP solves 76 instances to optimality that FFLBBD does not. Thus, our DIDP approach solves 67 more instances to optimality than FFLBBD while DIDP closes 76 open instances (i.e., DIDP failed to optimally solve 9 instances which FFLBBD did solve).

For data set A, B, and C, DIDP outperforms FFLBBD in terms of the average optimality gap (`Gap' column), average runtime (`Time' column), the number of instances with feasible solutions found (`Feas' column), and the number of instances with optimal solutions found and proved (`Opt' column). However, for the extra large data set D, FFLBBD is better, though DIDP shows some advantages in proving optimality and solution speed for $\alpha=0.75$ and $\alpha=1.00$.

\begin{table}[htbp]
  \centering
  \caption{Results of DIDP and FFLBBD: better results are in bold.}
  \label{sota_no_better}
  \setlength{\tabcolsep}{2.8mm}{}{
  \begin{tabular}{lllrrrrrrrr}
  \toprule[1.0pt]
  \multirow{2}{*}{Class} & \multirow{2}{*}{$\alpha$} & \multirow{2}{*}{\#} & \multicolumn{4}{c}{DIDP} & \multicolumn{4}{c}{FFLBBD} \\
  \cmidrule(lr){4-7} \cmidrule(lr){8-11} 
  & & & Gap & Time & Feas & Opt & Gap & Time & Feas & Opt \\
  \midrule[1.0pt]
  \multirow{4}{*}{A} 
  & 0.25  & 33 & \textbf{0.00} & \textbf{0.0007} & \textbf{33} & \textbf{33} & \textbf{0.00} & 0.12 & \textbf{33} & \textbf{33} \\
  & 0.50  & 33 & \textbf{0.00} & \textbf{0.0012} & \textbf{33} & \textbf{33} & \textbf{0.00} & 0.49 & \textbf{33} & \textbf{33} \\
  & 0.75  & 33 & \textbf{0.00} & \textbf{0.0022} & \textbf{33} & \textbf{33} & \textbf{0.00} & 0.64 & \textbf{33} & \textbf{33} \\
  & 1.00  & 33 & \textbf{0.00} & \textbf{0.0038} & \textbf{33} & \textbf{33} & \textbf{0.00} & 0.74 & \textbf{33} & \textbf{33} \\
  \midrule[1.0pt]
  \multirow{4}{*}{B} 
  & 0.25  & 35 & \textbf{0.00} & \textbf{3} & \textbf{35} & \textbf{35} & 0.02 & 28 & \textbf{35} & 33 \\
  & 0.50  & 35 & \textbf{0.00} & \textbf{4} & \textbf{35} & \textbf{35} & 0.01 & 42 & \textbf{35} & 34 \\
  & 0.75  & 35 & \textbf{0.00} & \textbf{7} & \textbf{35} & \textbf{35} & 0.10 & 51 & \textbf{35} & 32 \\
  & 1.00  & 35 & \textbf{0.00} & \textbf{8} & \textbf{35} & \textbf{35} & 0.11 & 54 & \textbf{33} & 33 \\
  \midrule[1.0pt]
  \multirow{4}{*}{C} 
  & 0.25  & 47 & \textbf{1.47} & \textbf{824} & \textbf{47} & \textbf{34} & 2.68 & 1399 & \textbf{47} & 13 \\
  & 0.50  & 47 & \textbf{3.07} & \textbf{922} & \textbf{47} & \textbf{28} & 5.71 & 1371 & \textbf{47} & 13 \\
  & 0.75  & 47 & \textbf{4.13} & \textbf{951} & \textbf{47} & \textbf{28} & 8.17 & 1350 & \textbf{47} & 14 \\
  & 1.00  & 47 & \textbf{4.99} & \textbf{919} & \textbf{47} & \textbf{27} & 10.39 & 1398 & \textbf{47} & 13 \\
  \midrule[1.0pt]
  \multirow{4}{*}{D} 
  & 0.25  & 82 & 11.25 & 1556 & \textbf{82} & 13 & \textbf{3.57} & \textbf{1440} & \textbf{82} & \textbf{18} \\
  & 0.50  & 82 & 14.03 & 1575 & \textbf{82} & 13 & \textbf{6.19} & \textbf{1551} & \textbf{82} & \textbf{15} \\
  & 0.75  & 82 & 17.14 & \textbf{1616} & \textbf{82} & \textbf{12} & \textbf{9.94} & 1620 & \textbf{82} & 11 \\
  & 1.00  & 82 & 18.77 & \textbf{1627} & \textbf{82} & \textbf{11} & \textbf{13.32} & 1639 & \textbf{82} & 9 \\
  \midrule[1.0pt]
  sum & & 788 & & & \textbf{788} & \textbf{438} & & & \textbf{788} & 371 \\
  \bottomrule[1.0pt]
  \end{tabular}}
\end{table}

\begin{table}[tp]
    \caption{Instances with different optimal objectives by using DIDP and FFLBBD.}
    \centering
    \setlength{\tabcolsep}{3.0mm}{}{
    \begin{tabular}{lllrrr}
        \toprule[1.0pt]
        Instance & Class & $\alpha$ & DIDP OptVal & FFLBBD UB & FFLBBD LB  \\
        \midrule[1.0pt]
        jackson\_c=7 & A & 0.75 & 9 & 8 & 8  \\
        jackson\_c=14 & A & 0.50 & 12 & 13 & 13  \\
        lutz1\_c=2357b & B & 0.25 & 2475 & 2472 & 2472  \\
        sawyer30\_c=54 & B & 0.75 & 58 & 57 & 57  \\
        hahn\_c=1806 & C & 0.25 & 2830 & 2848 & 2848  \\
        hahn\_c=4676 & C & 0.25 & 4847 & 4798 & 4798  \\
        \bottomrule[1.0pt]
    \end{tabular}}
    \label{table:diff}
\end{table}

However, our results indicate some problems with the results published in the online appendix of Zohali et al. \citeyearpar{zohali2022solving} where they provide the objective function values for each problem instance. For the small instance \texttt{jackson\_c=14.alb}, the best lower bound (LB) and upper bound (UB) on the cycle time according to Zohali et al. are both 13. Our DIDP, CP, and MIP models all find a solution with cycle time of 12. We have manually confirmed that the  DIDP solution is feasible, implying an error in the LB provided by Zohali et al. For the small instance \texttt{jackson\_c=7.alb}, the best LB and UB of the cycle time according to Zohali et al. are 8. Our DIDP, CP, and MIP models all find an optimal solution with cycle time being 9.  We do not have access to the detailed solution produced by FFLBBD for this instance and so are unable to analyze it further. There are six instances with inconsistent optimal objectives obtained by DIDP/CP/MIP and FFLBBD, as shown in Table \ref{table:diff}. In both cases where DIDP finds a better `optimal' solution, we have verified the solution feasibility. Thus, some of the results reported by Zohali et al. are incorrect.

\subsection{Results of the Local Improvement Algorithm for SUALBP-2}

Interstingly, as we can see in Figure \ref{time2} and \ref{pi2}, the performance of the local improvement algorithm is worse than our core DIDP model. To explain this observation, we use two metrics: the mean number and mean proportion of incumbent solutions that are improved by local improvement in each problem instance.

\begin{table}[tp]
    \caption{Mean number and proportion of CABS iterations that lead to real improvement.}
    \centering
    \setlength{\tabcolsep}{3.0mm}{}{
    \begin{tabular}{l|ccccc}
        \toprule[1.0pt]
        Instances & All & A & B & C & D  \\
        \midrule[1.0pt]
        Mean number & 0.897 & 0.167 & 0.386 & 0.686 & 1.530 \\
        Mean proportion & 0.081 & 0.037 & 0.053 & 0.071 & 0.117 \\
        \bottomrule[1.0pt]
    \end{tabular}}
    \label{table:licabs}
\end{table}

The results are shown in Table \ref{table:licabs}. Over all instances, local improvement finds a new incumbent 8.1\% of the times that it is called, which appears non-trivial. 
However, as local improvement solves the subproblems to optimality, it also means that in 91.9\% of the calls, our core model had already found the optimal sequence on the station with the maximum cycle time (given the task assignment).
As the local improvement is only called when an incumbent is found, on average it only finds 0.897 improving solutions per instance. Thus, local improvement does not bring performance gain with the runtime spent on it.
The mean number and the value of local improvement is increasing as the problem size gets larger, and, in fact, local improvement does find better solutions than our core model for some of the extra large instances.

\section{Conclusions}

In this paper, we study the assembly line balancing problem with sequence-dependent setup time (SUALBP), where tasks need to be assigned to stations and sequenced within each station, taking into account the processing time of tasks and sequence-dependent setup time between tasks. In type-1 SUALBP, an upper limit of station time is given and the objective is to minimize the number of active stations, while in type-2 SUALBP, the number of stations is fixed and the objective is to minimize the maximum station time.

We developed novel domain-independent dynamic programming (DIDP) and constraint programming (CP) optimization models. The DIDP models are formulated as state-based transition systems and solved with a DIDP solver using complete anytime beam search. We also proposed state-based dual bounds to accelerate the solving process of DIDP models. The CP models adopt a scheduling perspective, using optional interval variables and sequencing constraints.

We compared the performance of the proposed DIDP and CP models with the state-of-the-art MIP models on a diverse set of SUALBP instances from the literature. Experimental results show that the DIDP models are  significantly superior to the CP and MIP models in terms of proving optimality and finding high-quality solutions. The CP models prove optimality for slightly fewer instances than MIP but can find higher quality solutions than MIP on average.

For type-2 SUALBP, we compare the proposed DIDP model with the state-of-the-art exact algorithm, the full-featured logic-based Benders decomposition (FFLBBD) \citep{zohali2022solving}. DIDP outperforms FFLBBD in terms of proving optimality, especially for small, medium, and large instances. For extra large instances, DIDP is competitive with FFLBBD. In particular, DIDP proves optimality for 67 more instances than FFLBBD and closes 76 open instances, with 9 instances that are proved by FFLBBD not solved to optimality by DIDP. We have also detected errors in the results provided by Zohali et al. \citeyearpar{zohali2022solving}, raising concerns about the correctness of their approach.

We investigate a local improvement addition to the DIDP model for SUALBP-2, where a station-specific DIDP model seeks to improve new global incumbent solutions by improving the task sequence without changing task-to-station assignments. However, as the re-sequencing rarely improves the solution quality, the overall performance of the algorithm is worse than the core DIDP model.

This work represents an important contribution to exact methods for SUALBP and demonstrates the promising prospect of DIDP for solving complex planning and scheduling problems. Our main directions for future work are the continued exploration of DIDP performance across diverse combinatorial optimization problems as well as a study of solver behavior to develop better DIDP solution approaches.

\ACKNOWLEDGMENT{%
This work is supported by the Natural Sciences and Engineering Research Council of Canada.
}

\begin{APPENDICES}

    \section{Transitions of the DIDP Models} \label{appendix:A}
    
    An alternative perspective on the DIDP models arises from the set of transitions that define the recursive structure. We present the transitions for our two SUALBP models here.
    
    \subsection{The DIDP Model for SUALBP-1} \label{appendix:A.1}
    
    \begin{table}[tp]
        \caption{Summary of the DIDP model for SUALBP-1.}
        \label{table:trans1}
        \centering
        \setlength{\tabcolsep}{3.0mm}{}{
        \begin{tabular}{l|ccc}
            \toprule[1.0pt]
            State & Type & Objects & Preference \\
            \midrule[1.0pt]
            $U$ & set & tasks $V$ &  \\
            $\kappa$ & integer resource &  & less  \\
            $p$ & element & tasks $V$ &  \\
            $f$ & element & tasks $V$ &  \\
            $r$ & integer resource & & more  \\
            \midrule[1.0pt]
            \midrule[1.0pt]
            Target state & \multicolumn{3}{c}{$U=V, \kappa=0, p=d_{s}, f=d_{s}, r=0$} \\
            Base case & \multicolumn{3}{c}{$U=\emptyset \wedge f=d_{s}$} \\
            Dual bound & \multicolumn{3}{c}{$\max 
    \begin{cases}
        \ceil[\Big]{ \frac{\underline{\mu}_{f} +  \sum_{i \in U} (\underline{\tau}_{i} + t_{i}) - (\overline{m} - \kappa) \cdot (\max_{i \in U}{\underline{\tau}_{i}}) + \max(\underline{m} - \kappa, 0) \cdot (\min_{i \in U}{\underline{\mu}_{i}}) - r}{c} } \\
        \ceil[\Big]{ \frac{\underline{\mu}_{f} +  \sum_{i \in U} t_{i} - r}{c} }  \\
        \ceil{ \frac{\sum_{i\in U} t_{i} - r}{c} } \\
        \sum_{i \in U} w_{i}^{2} + \ceil{ \sum_{i\in U} w_{i}^{'2} - l^{2} } \\
        \ceil{ \sum_{i\in U} w_{i}^{3} - l^{3} } 
    \end{cases}$} \\
            \midrule[1.0pt]
            \midrule[1.0pt]
            Transition & Preconditions & Effects & Cost  \\
            \midrule[1.0pt]
            $\texttt{assign\_first}_{i}$ & 
            \makecell{$i \in V$, $f = d_{s}$, \\ $U \cap P_{i}^{*} = \emptyset$} & 
            \makecell{$U \rightarrow U \backslash \{i\}$, $r \rightarrow c - t_{i}$, \\ $p \rightarrow i$, $f \rightarrow i$, $\kappa \rightarrow \kappa+1$} & 
            1 \\
            \hline
            $\texttt{assign\_next}_{i}$ & 
            \makecell{$i \in U$, $U \cap P_{i}^{*} = \emptyset$, \\ $f \neq d_{s}$, $t_{i} + \tau_{p i} \leq r$} & 
            \makecell{$U \rightarrow U \backslash \{i\}$, $p \rightarrow i$, \\ $r \rightarrow c - t_{i} - \tau_{p i}$} & 
            0 \\
            \hline
            $\texttt{close\_station}$ & 
            \makecell{$\{i \notin U$ \text{ or } $t_{i} + \tau_{p i} + \mu_{i f} > r $ \\ $\text{ or } U \cap P_{i}^{*} \neq \emptyset, \forall i \in V\}$, \\ $f \neq d_{s}$, $\mu_{p f} \leq r$} & 
            $r \rightarrow 0$, $p \rightarrow d_{s}$, $f \rightarrow d_{s}$ & 
            0 \\
            \bottomrule[1.0pt]
        \end{tabular}}
    \end{table}
    
    The transitions of the DIDP model for SUALBP-1 are presented in Table \ref{table:trans1}. The first type of transition is opening a new station and assigning task $i$ as its first task. The second type of transition is assigning task $i$ as the next task on the current station that already has at least one task assigned. The forward setup time to task $i$ is handled with this transition. The third type of transition is closing the current station. The backward setup time from task $p$ to task $f$ is handled with this transition. The precondition of the final transition ensures that there is no task in $U$ that can be assigned as the last task of the current station. However, it does not mean we cannot assign more tasks to the current station.
    It is possible that $\exists i,j \in U$ such that $t_{j} + \tau_{i j} + \mu_{j f} < \mu_{i f}$ and hence $t_{i} + \tau_{p i} + t_{j} + \tau_{i j} + \mu_{j f} < t_{i} + \tau_{p i} + \mu_{i f}$. As a result, we can assign tasks $i$ and $j$ to the current state without exceeding the remaining time $r$. In order to incorporate this possibility, we allow transition $\texttt{assign\_next}_{i}$ to consider any task $i$ as long as $i \in U$, $U \cap P_{i}^{*} = \emptyset$, and $t_{i} + \tau_{p i} \leq r$, where backup setup times are not involved. At the same time, we endow the transition \texttt{close\_station} with the freedom to close the current station if no task in the unscheduled task set can be selected as the next task for the current station. If no transition $\texttt{assign\_next}_{i}$ can occur at the current state and $\mu_{p f} > r$, the current state is a dead-end.

    \begin{table}[htbp]
        \caption{Summary of the DIDP model for SUALBP-2.}
        \label{table:trans2}
        \centering
        \setlength{\tabcolsep}{3.0mm}{}{
        \begin{tabular}{l|ccc}
            \toprule[1.0pt]
            State & Type & Objects & Preference  \\
            \midrule[1.0pt]
            $U$ & set & tasks $V$ &  \\
            $\kappa$ & integer &  &  \\
            $p$ & element & tasks $V$ &  \\
            $f$ & element & tasks $V$ &  \\
            $t^{c}$ & integer resource & & less  \\
            $c$ & integer resource & & less  \\
            \midrule[1.0pt]
            \midrule[1.0pt]
            Target state & \multicolumn{3}{c}{$U=V, \kappa=0, p=d_{s}, f=d_{s}, t^{c}=0, c=0$} \\
            Base case & \multicolumn{3}{c}{$U=\emptyset \wedge f=d_{s}$} \\
            Dual bound & \multicolumn{3}{c}{$\max
    \begin{cases}
        \ceil[\Big]{ \frac{ \sum_{i \in U} (\underline{\tau}_{i} + t_{i}) + t^{c} + (m - \kappa) \cdot (\min_{i \in U}{\underline{\mu}_{i}} - \max_{i \in U}{\underline{\tau}_{i}}) + \underline{\mu}_{f}}{\min (m, m - \kappa + 1)} - c}  \\
        \ceil[\Big]{ \frac{ \sum_{i \in U} t_{i} + t^{c} }{\min (m, m - \kappa + 1)} - c}
    \end{cases}$} \\
            \midrule[1.0pt]
            \midrule[1.0pt]
            Name & Preconditions & Effects & Cost  \\
            \midrule[1.0pt]
            $\texttt{assign\_first}_{i}$ & 
            \makecell{$i \in U$, $f = d_{s}$, \\ $U \cap P_{i}^{*} = \emptyset$, $\kappa < m$} & 
            \makecell{$U \rightarrow U \backslash \{i\}$, $t^{c} \rightarrow t_{i}$, \\ $p \rightarrow i$, $f \rightarrow i$, $\kappa \rightarrow \kappa + 1$, \\ $c \rightarrow \max{(c, t_{i})}$} & 
            $\max{(0, t_{i} - \textcolor{purple}{c})}$ \\
            \hline
            $\texttt{assign\_next}_{i}$ & 
            \makecell{$i \in U$, $f \neq d_{s}$, \\ $U \cap P_{i}^{*} = \emptyset$} & 
            \makecell{$U \rightarrow U \backslash \{i\}$, $p \rightarrow i$, \\ $t^{c} \rightarrow t^{c} + t_{i} + \tau_{p i}$, \\ $c \rightarrow \max{(c, t^{c} + t_{i} + \tau_{p i})}$} & 
            \makecell{$\max(0, t^{c} +$ \\ $t_{i} + \tau_{p i} - \textcolor{purple}{c})$} \\
            \hline
            $\texttt{close\_station}$ & 
            \makecell{$f \neq d_{s}$} & 
            \makecell{$t^{c} \rightarrow t^{c} + \mu_{p f}$, $p \rightarrow d_{s}$, \\ $f \rightarrow d_{s}$, $c \rightarrow \max{(c, t^{c} + \mu_{p f})}$} & 
            \makecell{$\max(0, t^{c} +$ \\ $\mu_{p f} - \textcolor{purple}{c})$} \\
            \bottomrule[1.0pt]
        \end{tabular}}
    \end{table}

    \subsection{The DIDP Model for SUALBP-2}  \label{appendix:A.2}

    The transitions in this model are similar to the DIDP model of SUALBP-1. In particular, we also use $\texttt{assign\_first}_{i}$, $\texttt{assign\_next}_{i}$, and $\texttt{close\_station}$. However, the state variables and costs are updated differently to account for the different cost function. The transitions are summarized in Table \ref{table:trans2}.

    \subsection{The DIDP Model for the sequencing subproblem of the local improvement algorithm} \label{appendix:A.3}

    In the sequencing subproblem of the local improvement algorithm, all tasks assigned to the same station are re-sequenced to minimize the total station time. The DIDP model for the problem is summarized in Table \ref{table:trans3}. $V_{k}$ represents the set of the tasks assigned to station $k$, $d_{k}$ represents the dummy task of station $k$, $P_{i}^{k}$ is the set of all the predecessors at station $k$ of task $i$, $\underline{\tau}_{i}^{k}$ is the smallest forward setup time from any task at station $k$ to task $i$, and $\underline{\mu}_{f}^{k}$ is the smallest backward setup time from any task at station $k$ to task $f$. Note that the processing time of tasks are excluded from the DIDP model as they are constant after the task assignment is fixed.
    
    There are three state variables. $U$ is a set variable representing unscheduled tasks \textcolor{purple}{and $U=V_{k}$} in the target state. $p$ is an element variable representing the previous task at station $k$ and $p=d_{k}$ in the target state. $f$ is an element variable representing the first task at station $k$ and $f=d_{k}$ in the target state. There are three types of transitions. The first type is assigning task $i$ as the first task at station $k$. The second type is assigning task $i$ as the next task at station $k$. This transition also handles the forward setup times. The third type is closing the station and it deals with the backward setup time.
    
    \noindent \emph{Base case.} The base case of the DIDP model is: $U = \emptyset \wedge p=d_{s}$. Note that $p=d_{s}$ is necessary since backward setup time has to be included in the base case, with the help of transition $\texttt{close\_station}$.

    \begin{table}[tp]
        \caption{Summary of the DIDP model for the sequencing subproblem.}
        \label{table:trans3}
        \centering
        \setlength{\tabcolsep}{3.0mm}{}{
        \begin{tabular}{l|ccc}
            \toprule[1.0pt]
            State & Type & Objects & Preference  \\
            \midrule[1.0pt]
            $U$ & set & tasks $V_{k}$ &  \\
            $p$ & element & tasks $V_{k}$ &  \\
            $f$ & element & tasks $V_{k}$ &  \\
            \midrule[1.0pt]
            \midrule[1.0pt]
            Target state & \multicolumn{3}{c}{$U=V, p=d_{k}, f=d_{k}$} \\
            Base case & \multicolumn{3}{c}{$U=\emptyset \wedge p=d_{k}$} \\
            Dual bound & \multicolumn{3}{c}{$\begin{cases}
                \sum_{i \in U} \underline{\tau}_{i}^{k} + \underline{\mu}_{f}^{k}, \text{ if } p \neq d_{k} \\
                0, \qquad \qquad \quad \ \text{ if } p = d_{k}
            \end{cases}$} \\
            \midrule[1.0pt]
            \midrule[1.0pt]
            Name & Preconditions & Effects & Cost  \\
            \midrule[1.0pt]
            $\texttt{assign\_first}_{i}$ & 
            \makecell{$i \in U$, $f = d_{k}$, $U \cap P_{i}^{k} = \emptyset$} & 
            \makecell{$U \rightarrow U \backslash \{i\}$, $p \rightarrow i$, $f \rightarrow i$} & 
            $0$ \\
            \hline
            $\texttt{assign\_next}_{i}$ & 
            \makecell{$i \in U$, $f \neq d_{k}$, $U \cap P_{i}^{k} = \emptyset$} & 
            \makecell{$U \rightarrow U \backslash \{i\}$, $p \rightarrow i$} & 
            \makecell{$\tau_{p i}$} \\
            \hline
            $\texttt{close\_station}$ & 
            \makecell{$U = \emptyset, f \neq d_{k}$} & 
            \makecell{$p \rightarrow d_{k}$, $f \rightarrow d_{k}$} & 
            \makecell{$\mu_{p f}$} \\
            \bottomrule[1.0pt]
        \end{tabular}}
    \end{table}

    \noindent \emph{Recursive function.} We use $\mathcal{V}(U, p, f)$ to represent the cost of a state. Let $U_{1} = \{j \in U | U \cap P_{j}^{k} = \emptyset\}$. The recursive function of the DIDP model is as follows:
    \begin{subequations} 
    \begin{align}
    &\texttt{compute} \ \mathcal{V}(V, d_{k}, d_{k}) \label{9a} &\\  
    & \mathcal{V}(U, p, f) = \label{9b} \\
    & \qquad \begin{cases}
        0 & \text{if $U=\emptyset, p=d_{k}$,}  \notag \qquad \qquad \ \ \text{(i)} \\
        \min_{j \in U_{1}} \mathcal{V}(U\backslash\{j\}, j, j) & \text{if $U_{1} \neq \emptyset, f=d_{k}$,} \qquad \qquad \ \text{(ii)} \\
        \min_{j \in U_{1}} (\tau_{p j} + \mathcal{V}(U\backslash\{j\}, j, f)) & \text{if $U_{1} \neq \emptyset, f \neq d_{k}$,} \qquad \qquad \ \text{(iii)} \\
        \min(\mu_{p f} + \mathcal{V}(U, d_{k}, d_{k})) & \text{if $U = \emptyset, f \neq d_{k}$,} \qquad \qquad \ \ \text{(iv)} \\
        \end{cases} \\ 
    & U_{1} = \{j \in U | U \cap P_{j}^{k} = \emptyset\}, \notag \\
    & \mathcal{V}(U, p, f) \geq 
    \label{9c} 
    \begin{cases}
        \sum_{i \in U} \underline{\tau}_{i}^{k} + \underline{\mu}_{f}^{k}, \text{ if } p \neq d_{k} \\
        0, \qquad \qquad \quad \ \text{ if } p = d_{k}
    \end{cases}
    \end{align}
    \end{subequations} 
    The term (\ref{9a}) is to compute the objective of the target state. Equation (\ref{9b}) is the main recursion of the DIDP model. Specifically, (20b-i) handles the base cases, while (20b-ii) refers to assigning the first task to the current station. (20b-iii) corresponds to assigning the next task after other tasks have been assigned. (20b-iv) deals with closing the current station. Term (\ref{9c}) formulates the state-based dual bounds. It uses $\underline{\tau}_{i}^{k}$ to underestimate the forward setup time to task $i \in U$ and $\underline{\mu}_{f}^{k}$ to underestimate the backward setup time to task $f$.

    \end{APPENDICES}

%
%
%


\bibliographystyle{informs2014} 
\bibliography{JOC-template.bib} 


\end{document}